\documentclass[a4wide]{article}
\pdfoutput=1 
\usepackage[english]{babel}
\usepackage{amsmath}
\usepackage{amssymb}
\usepackage{theorem}
\usepackage{multibbl}

\usepackage{afterpage}
\usepackage[retainorgcmds]{IEEEtrantools}
\usepackage{graphicx}
\usepackage{subfigure}
\graphicspath{{figss/}}
\usepackage{color}
\usepackage{pstricks,pst-all,pstricks-add,pst-plot,pst-eucl}
\usepackage{enumerate}

\theorembodyfont{\it}
\newtheorem{theorem}{Theorem}

\newtheorem{lemma}{Lemma}

\theorembodyfont{\rm}

\newtheorem{remark}{Remark}

\makeatletter
\renewcommand{\fnum@figure}{\small\textbf{\figurename~\thefigure}}
\makeatother

\definecolor{brown}{rgb}{0.64,0.16,0.16}
\definecolor{ForestGreen}{rgb}{0.13,0.54,0.13}
\definecolor{purple}{rgb}{0.62,0.12,0.94}
\definecolor{DodgerBlue}{rgb}{0.11,0.56,0.98}
\definecolor{RoyalBlue}{rgb}{0.25,0.41,0.88}
\definecolor{B}{rgb}{0,0,1}
\definecolor{G}{rgb}{0,0.502,0}
\definecolor{R}{rgb}{1,0,0}

\def \RR {\mathbb{R}}

\title{A Novel Phase Portrait to Understand Neuronal Excitability}
\author{Alessio Franci$^{1,*}$, Guillaume Drion$^{2,3,*}$, Vincent Seutin$^{2}$ \& Rodolphe Sepulchre$^{3}$\\
\small{$^1$L2S, University Paris Sud 11, Gif-sur-Yvette, France.}\\
\small{$^2$Laboratory of Pharmacology and GIGA Neurosciences, University of Li\`ege, Li\`ege, Belgium.}\\
\small{$^3$Department of Electrical Engineering and Computer Science \& GIGA-research, University of Li\`ege, Li\`ege, Belgium.}\\
\small{*These authors contributed equally to this work.}}
\date{}

\begin{document}

\maketitle

\begin{abstract}
Fifty years ago, Fitzugh introduced a phase portrait that became famous for a twofold reason: it captured in a physiological way the qualitative behavior of Hodgkin-Huxley model and it revealed the power of simple dynamical models to unfold complex firing patterns. To date, in spite of the enormous progresses in qualitative and quantitative
neural modeling, this phase portrait has remained the core picture of neuronal excitability. Yet, a major difference between the neurophysiology of 1961 and of 2011 is the recognition of the prominent role of calcium channels in firing mechanisms. We show that including this extra current  in Hodgkin-Huxley dynamics leads to a revision of Fitzugh-Nagumo phase portrait that affects in a fundamental way the reduced modeling of neural excitability. The revisited model considerably enlarges the modeling power of the original one. In particular, it captures essential electrophysiological signatures that otherwise require non-physiological alteration or considerable complexication of the classical model. As a basic illustration, the new model is shown to highlight a core dynamical mechanism by which the calcium conductance controls the two distinct firing modes of thalamocortical neurons.
\end{abstract}

\newbibliography{main}

\section{Introduction}
Rooted in the seminal work of Hodgkin and Huxley \cite{main}{HODHUX}, conductance-based models have become a central paradigm to describe the electrical behavior of neurons. These models
combine a number of advantages, including physiological interpretability (parameters have a precise
experimental meaning) and modularity (additional ionic currents and/or spatial effects are easily
incorporated using the interconnection laws of electrical circuits \cite{main}{halnes2011multi,canavier2006increase}).
Not surprisingly, the gain in quantitative description is achieved at the expense of mathematical
complexity. The dimension of detailed quantitative models makes them mathematically intractable for analysis and numerically intractable for the simulation of large neuronal populations. For this reason, reduced modeling of
conductance-based models has proven  an indispensable complement to quantitative
modeling. In particular, the FitzHugh-Nagumo model \cite{main}{fitzhugh61}, a two-dimensional reduction of Hodgkin-Huxley model, has played an essential role over the years to explain the mechanisms of neuronal excitability (see e.g. \cite{main}{Rinzel:1989:ANE:94605.94613, ERTE10} for an excellent introduction and further references). More recently, Izhikevich has enriched the value of reduced-models by providing the Fitzugh-Nagumo model with a reset mechanism \cite{main}{IZHIKEVICH03} that captures the fast (almost discontinuous) behavior of spiking neurons. Such models are used to reproduce the qualitative \cite{main}{IZHIKEVICH2010,IZHIKEVICH2007} and quantitative \cite{main}{pospischil2011comparison,richert2011efficient} behavior of a large family of neuron types.  Notably, their computational economy makes them good candidates for large-scale simulations of neuronal populations \cite{main}{izhikevich2008large}.
 
The Hodgkin-Huxley model and all reduced models derived from it \cite{main}{fitzhugh61,IZHIKEVICH2007} focus on sodium and potassium currents, as the main players in the generation of action potentials: sodium is a fast depolarizing current, while potassium is slower and hyperpolarizing. 
Initally motivated by reduced modeling of dopaminergic neurons in which calcium currents are essential to the firing mechanisms \cite{main}{DRMASESE11}, the present paper mimicks the classical reduction of the Hodgkin-Huxley model augmented with an additional calcium current. The calcium current is a distinct player because it is depolarizing, as the sodium current,  but acts on the slower timescale of the potassium current.
 
To our surprise, the inclusion of calcium currents in the HH model before its planar reduction leads to a novel phase portrait that has been disregarded to date. Mimicking earlier classical work, we perform a normal form reduction of the global HH reduced planar model. The mathematical normal form reduction is fundamentally different in the classical and new phase portrait because it involves a different bifurcation. The classical fold bifurcation is replaced by a transcritical bifurcation.

The results of these mathematical analysis lead to a novel simple model that further enriches the modeling power of the popular hybrid model of Izhikevich. A single parameter in the new model controls the neuron calcium conductance. In low calcium conductance mode, the model captures the standard behavior of earlier models. But in high calcium conductance mode, the same model captures the electrophysiological signature of neurons with a high density of calcium channels, in agreement with many experimental observations. For this reason, the novel reduced model is particularly relevant to understand the firing mechanisms of neurons that switch from a low calcium-conductance mode to a high calcium-conductance mode. Because thalamocortical (TC) neurons provide a prominent example of such neurons, they are chosen as the main illustration of the present paper, the benefits of of which should extend to a much broader class of neurons.

\section{Planar reduction of Hodgkin-Huxley model revisited in the light of calcium channels}
\label{SEC: proposed model and planar reduction}

Calcium channels participate in the spiking pattern by providing, together with sodium channels, a source of depolarizing currents. In contrast to sodium channels whose gating kinetics are fast, calcium channels activate on a slower time-scale, similar to that of potassium channels \cite{main}{hille1991ionic}. As a consequence, they oppose the hyperpolarizing effect of potassium current activation, resulting in bidirectional modulation capabilities of the post-spike refractory period. We model this important physiological feature by considering the HH model \cite{main}{HODHUX} with an additional voltage-gated (non-inactivating) calcium current $I_{Ca}$ and a DC-current $I_{pump}$ that accounts for hyperpolarizing calcium pump currents. The resulting model is similar to HH model when the conductance of the calcium currents is low, but becomes strikingly different when the calcium conductance is high. Note that the inactivation of calcium channels is not included in the HH dynamics because it is known to take place in a much slower time scale \cite{main}{wang1991model}. The inactivation will typically be accounted for by a slower adaptation of the calcium conductance, see Section \ref{SEC: TC hybrid modeling}.

Figure 1{\bf a} illustrates the spiking behavior induced by the action of an external square current $I_{app}$ in the two different modes. As compared to the original HH model (Fig. 1{\bf a} left), the presence of the calcium current is characterized by a triple electrophysiological signature (see Fig. 1{\bf a} right):
\begin{itemize}
\item \underline{spike latency}: the spike train (burst) is delayed with respect to the onset of the stimulation
\item \underline{plateau oscillations}: the spike train oscillations occur at a more depolarized voltage than the hyperpolarized state
\item \underline{after-depolarization potential} (ADP): the burst terminates with a small depolarization
\end{itemize}
This electrophysiological signature is typical of neurons with sufficiently strong calcium currents. See for instance: spike latency \cite{main}{Rekling01111997,Molineux23112005}, plateau oscillations \cite{main}{Beurrier15011999}, ADPs \cite{main}{Azouz01041996,TJP:TJP2775}. However, the mechanisms by which these behaviors occur have never been analyzed using reduced planar models to date.

\begin{figure}
\centering
\includegraphics[width=0.7\linewidth]{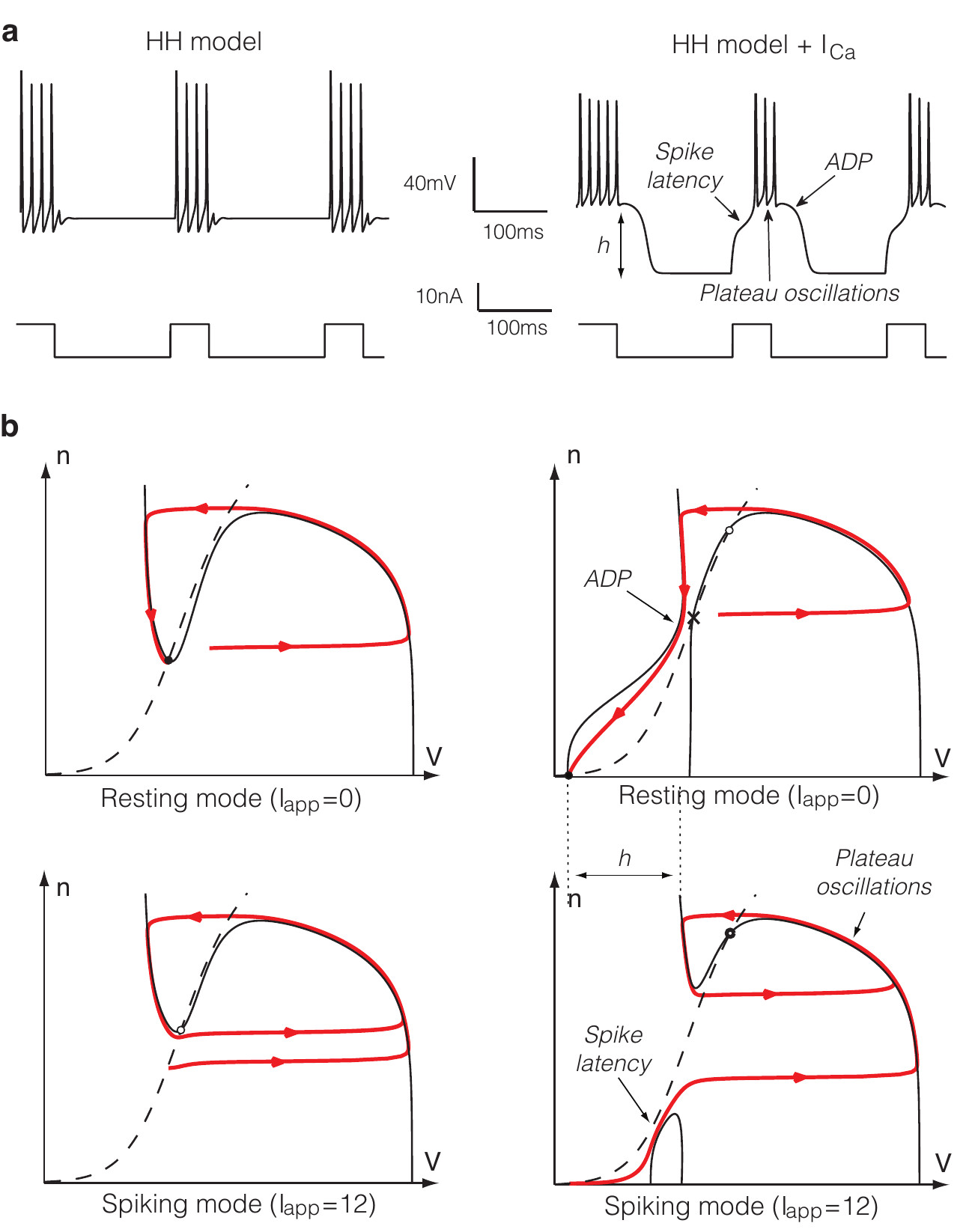}
\caption{{\bf Step responses of the HH model without (left) and with a calcium current (right).} ({\bf a}) Time-evolution of the applied excitatory current (bottom) and of the corresponding membrane potential (top) in HH model (the reduced model leads to almost the same behavior (Supplementary Fig. S1)). ({\bf b}) Phase portraits of the reduced Hodgkin-Huxley model in resting (top) and spiking states (bottom). The $V$- and $n$-nullclines are drawn as a full and a dashed line, respectively. Trajectories are drawn as solid oriented red lines. Black circles denote stable fixed points, white circles unstable fixed points, and cross saddle points. The presence of calcium channels strongly affects the phase-portrait and the corresponding  electrophysiological time-response of the neuron to excitatory inputs.}\label{FIG:1}
\end{figure}

Following the standard reduction of HH model \cite{main}{fitzhugh61}, we concentrate on the voltage variable $V$ (that accounts for the membrane potential) and on a recovery variable $n$ (that accounts for the overall gating of the ion channels) as key variables governing excitability (see methods). The phase-portrait of the reduced HH model is shown in Figure 1{\bf b} (left). This phase portrait and the associated reduced dynamics are well studied in the literature (see \cite{main}{fitzhugh61} for the FitzHugh paper, and \cite{main}{ERTE10,IZHIKEVICH2007} for a recent discussion and more references). We recall them for comparison purposes only. The resting state is a stable focus, which lies near the minimum of the familiar N-shaped $V$-nullcline. When the stimulation is turned on (spiking mode), this fixed point loses stability in a subcritical Andronov-Hopf bifurcation (see the numerical bifurcation analysis of Supplementary Section \ref{SEC: bif analysis}), and the trajectory rapidly converges to the periodic spiking limit cycle attractor. As the stimulation is turned off (resting mode), the resting state recovers its global attractivity via a saddle-node of limit cycles (the unstable one being born in the subcritical Hopf bifurcation), and the burst terminates with small subthreshold oscillations (cf. Fig. 1{\bf a} left).

In the presence of the calcium current, the phase-portrait changes drastically, as shown in Figure 1{\bf b} (right). In the resting mode, the hyperpolarized state is a stable node lying on the far left of the phase-plane. The $V$-nullcline exhibits a ``hourglass'' shape. Its left branch is attractive and guides the relaxation toward the resting state after a single spike generation. The sign of $\dot V$ changes from positive to negative approximately at the funnel of the hourglass, corresponding to the ADP apex. The right branch is repulsive and its two intersections with the $n$-nullcline are a saddle and an unstable focus.

When the stimulation is turned on, the $V$-nullcline breaks down in an upper and a lower branch. The upper branch exhibits the familiar N-shape and contains an unstable focus surrounded by a stable limit cycle, very much as in the reduced Hodgkin-Huxley model. In contrast, the lower branch of the $V$-nullcline, which is not physiological without the calcium currents, comes into play. While converging toward the spiking limit cycle attractor from the initial resting state, the trajectory must travel between the two nullclines where the vector field has smaller amplitude. As a consequence, the first spike is fired with a latency with respect to the onset of the stimulation, as observed in Figures 1{\bf a} (right) in the presence of the calcium current (see also Supplementary Fig. S\ref{FIG: red comparison}). In addition, a comparison of the relative position of the resting state and the spiking limit cycle in Figure 1{\bf b} (right) explains the presence of plateau oscillations. As the stimulation is turned off the spiking limit cycle disappears in a saddle-homoclinic bifurcation (see Supplementary Sections \ref{SEC: bif analysis} and \ref{SEC: trans and homo bifurcations}), and the resting state recovers its attractivity.

The presence of the lower branch of the $V$-nullcline has a physiological interpretation. In the reduced HH model, the gating variable $n$ accounts for the activation of potassium channels and the inactivation of sodium channels. Their synergy results in a total ionic current that is monotonically increasing with $n$ for a fixed value of $V$ (Supplementary Fig. S\ref{EQ: V-null phys interpret} left). In this situation, at most one value of $n$ solves the equation $\dot V=0$ and there can be only one branch for the voltage nullcline. In contrast, when calcium channels are present, the reduced gating variable must capture two antagonistic effects. As a result, the total ionic current is decreasing for low $n$ (the gating variable is excitatory), and increasing for large $n$ (the gating variable recovers its inhibitory nature) (Supplementary Fig. S\ref{EQ: V-null phys interpret} right). In this situation, two distinct values of $n$ solve the equation $\dot V=0$, which explains physiologically the second branch of the $V$-nullcline. To summarize, the lower branch of the voltage nullcline accounts for the existence of an excitatory effect of $n$, which is brought by calcium channel activation.

\section{The central ruler of excitability is a transcritical bifurcation, not a fold one}

The power of mathematical analysis of the reduced planar model (\ref{EQ: redHHCA}) is fully revealed by introducing two further simplifications.

\begin{itemize}
\item Time-scale separation: we exploit that the voltage dynamics are much faster than the recovery dynamics by assuming a small ratio $\dot n = O(\varepsilon) \dot V$ (the approximation holds away from the voltage nullcline) and by studying the singular limit $\varepsilon = 0$.
\item Transcritical singularity: by comparing the shape of the voltage nullcline in Fig. 1(b) ($I = 0$) and Fig. 1(d) ($I = 12$), one deduces from a continuity argument that a critical value $0 < I_c < 12$ exists at which the two branches of the voltage nullcline intersect.
\end{itemize}

The critical current $I_c$ depends on $\epsilon$. In the singular limit ($\epsilon = 0$) and for the corresponding critical current $I_c = I_c(0)$, one obtains the highly degenerate phase portrait in Figure 2{\bf a} (center). This particular phase portrait contains a transcritical bifurcation (red circle) which is the key ruler of excitability. This is because,  as illustrated in Figure 2{\bf b} for $\epsilon>0$,  the convergence of solutions  either to the resting point ($I<I_c$) or to the spiking limit cycle ($I>I_c$) is fully determined by the stable $W_s$ and unstable $W_u$ manifolds of the saddle point. In the singular limit shown in Figure 2{\bf a}, these hyberbolic objects degenerate to a critical manifold that coincides with the voltage nullcline near the transcritical bifurcation. It is in that sense that the X-shape of the voltage nullcline completely organizes the excitability, i.e. the transition from resting state to limit cycle.

\begin{figure}
\centering
\includegraphics[width=0.8\linewidth]{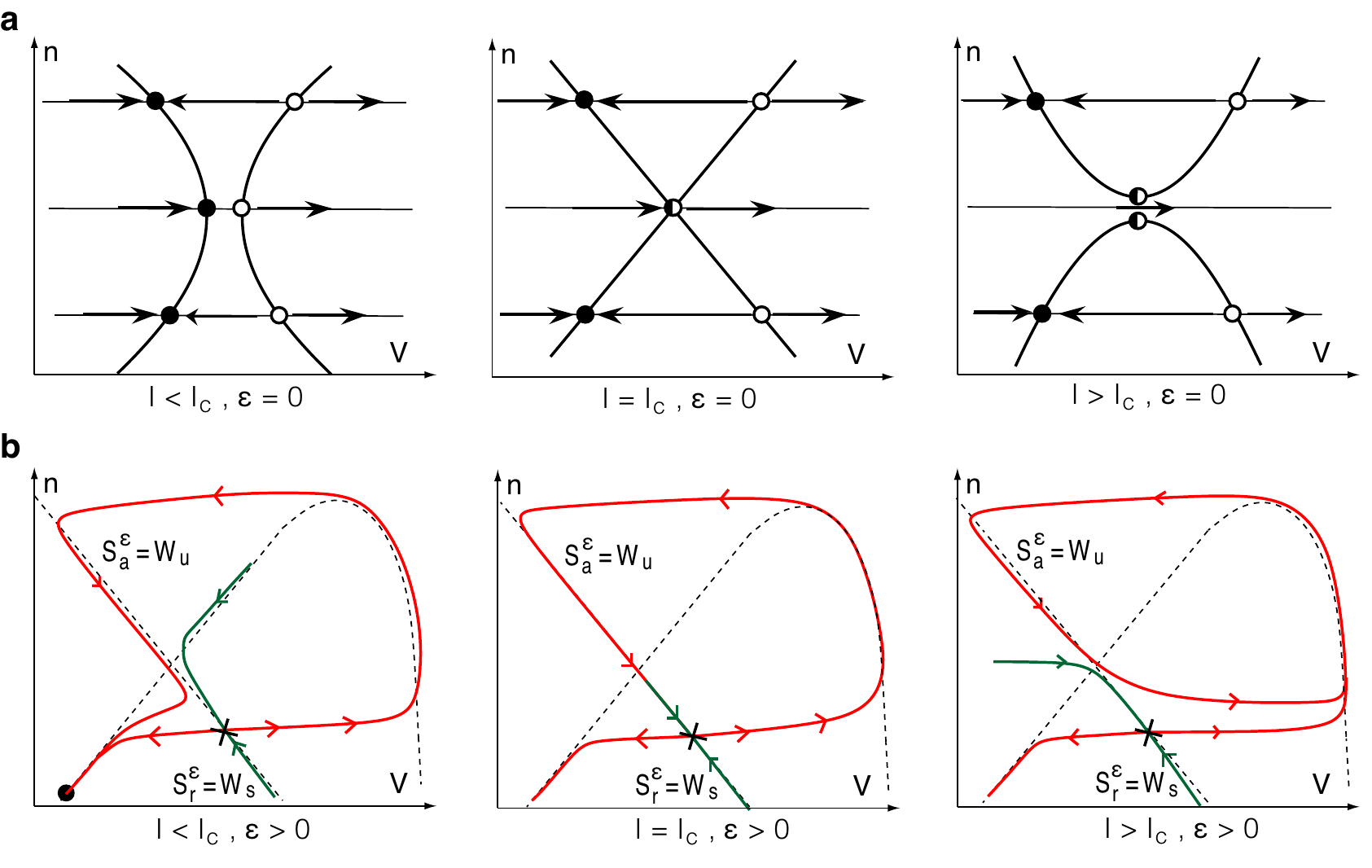}
\caption{{\bf Transcritical bifurcation as the main ruler of neuronal excitability.} ({\bf a}) Cartoon of the V-nullcline transition through a singularly perturbed transcritical bifurcation. Black circles denote stable fixed points, white circles unstable fixed points. ({\bf b}) Continuation of the stable $W_s$ (in green) and the unstable $W_u$ (in red) manifolds of the saddle away from the singular limit ({\it i.e.} $\epsilon>0$). They dictate the transition from the resting state ($I < I_c$) to the the spiking limit cycle ($I > I_c$) via a saddle-homoclinic bifurcation ($I=I_c$).}\label{FIG: trans}
\end{figure}

The persistence of the manifold $W_s$ and $W_u$ away from the singular limit is proven by geometric singular perturbation (Supplementary Section \ref{SEC: trans and homo bifurcations}). The same analysis also establishes a normal form behavior in the neighborhood of the transcritical bifurcation: in a system of local coordinates centered at the bifurcation, the voltage dynamics take the simple form
\begin{eqnarray}\label{EQ: trans normal form temp}
\dot v&=&v^2-w^2+i+h.o.t.\nonumber
\end{eqnarray}
where $i$ is a re-scaled input current and with $h.o.t.$ referring to higher order terms in $v,w,\varepsilon$.

It should be emphasized that it is the same perturbation analysis that leads to the classical view of the Hodgkin-Huxley reduced dynamics: the transition from Figure 1{\bf b} left ($I = 0$) to Figure 1{\bf b} left ($I = 12$) involves a fold bifurcation that governs the excitability with a fold normal form
\begin{eqnarray}
\dot v&=&v^2-w+i+h.o.t.\nonumber
\end{eqnarray}

\begin{figure}
\centering
\includegraphics[width=0.8\linewidth]{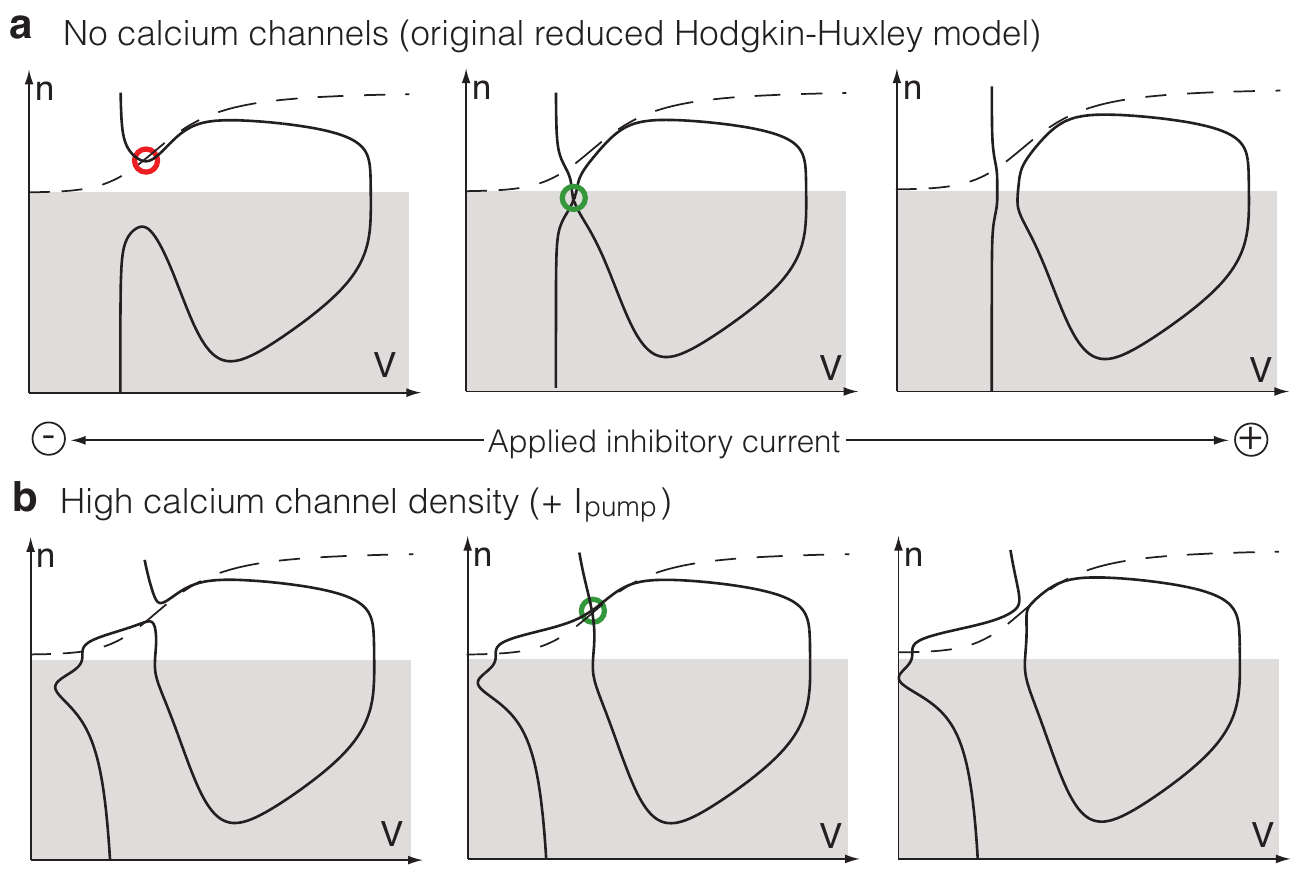}
\caption{{\bf Unfolding of the transcritical bifurcation in the global reduced Hodgkin-Huxley phase portrait.} ({\bf a} and {\bf b}) Phase portraits or the original reduced HH model without ({\bf a}) and with a calcium current ({\bf b}). A constant inhibitory current of increasing amplitude (from left to right) is applied to the model. The transcritical bifurcation is non physiological in the classical reduced HH model ({\bf a}) but plays an important physiological role in the presence of calcium channels ({\bf b}).}\label{FIG: global phase portrait}
\end{figure}

It is of interest to realize that the addition of the calcium current in the HH model unmasks a global view of its phase portrait that has been disregarded to date for its lack of physiological relevance. Figure \ref{FIG: global phase portrait} (top) shows the phase portrait of the classical reduced HH model for three different values of the hyperpolarizing current, revealing the transcritical singularity for the middle current value. The unshaded part of the first plot (and only this part of the plot) is familiar to most neuroscientists since the work of FitzHugh. Likewise, the conceptual sketch of the transcritical bifurcation will be familiar to all readers of basic textbooks in bifurcation analysis. For instance, the sketch is found in \cite{main}{SEYDEL94} as a prototypical example of non-generic bifurcation. It is symptomatic that this particular example is described at length but not connected to any concrete model in a texbook that puts much emphasis on the relevance of bifurcation analysis in neurodynamics applications. As shown in Fig \ref{FIG: global phase portrait} (bottom), the missing connection is brought to life by calcium channels. Their particular kinetics renders the transcritical bifurcation of HH model physiological in the the presence of a high-conductance calcium current.

\section{Transcritical hybrid modeling of neurons}
\label{SEC: new reduced hybrid model}

The singular limit of planar reduced models reveals that the excitability properties of spiking neurons are essentially determined by a local normal form of bifurcation of the resting equilibrium. This property is at the core of mathematical analysis of neuronal excitability (see \cite{main}{ERTE10,IZHIKEVICH2007} and the rich literature therein).

In recent work, Izhikevich showed that, for computational purposes, the combination of the local normal form dynamics with a hybrid reset mechanism, mimicking the fast (almost discontinuous) spike down-stroke, is able to reproduce the behavior of a large family of neurons with a high degree of fidelity  \cite{main}{IZHIKEVICH03,IZHIKEVICH2010}. Mimicking Izhikevich approach, we simplify the planar dynamics into the hybrid model:
\begin{IEEEeqnarray}{rCl}\label{EQ: IZHw2 model}
\dot v=v^2-w^2+I &\quad\quad\quad\mbox{if } v \geq v_{th}, \mbox{then}\IEEEyessubnumber\\
\dot w=\epsilon(av-w+w_0) &\quad\quad\quad v\leftarrow c, w\leftarrow d\IEEEyessubnumber
\end{IEEEeqnarray}

The proposed transcritical hybrid model is highly reminiscent of the hybrid model of Izhikevich, but it consideraly enlarges its modeling power by including two features of importance:

\begin{itemize}
\item the transcritical normal form $\dot{v} = v^2 - w^2 + I$ replaces the fold normal form $\dot{v} = v^2 - w + I$, in accordance with the normal form analysis of Section 3.
\item the new parameter $w_0$ determines whether the intersection of the voltage and recovery nullclines will take place above ($w_0 > 0$) or below ($w_0 < 0$) the transcritical singularity.
\end{itemize}

The parameter $w_0$ is a direct image of the calcium conductance: for small calcium conductances, the recovery variable nullcline only intersects the upper branch of the voltage nullcline (Supplementary Fig. S\ref{FIG: Na to Ca}) ; likewise in the hybrid model when $w_0 > 0$. For high calcium conductance, the recovery variable nullcline intersects the lower branch of the voltage nullcline (Supplementary Fig. S\ref{FIG: Na to Ca}) ; likewise in the hybrid model when $w_0 < 0$. 

Fig. 4 summarizes the four different phase portraits that derive from the transcritical hybrid model for different values of $I$ and $w_0$. For $w_0 > 0$, the model captures the classical view of the reduced HH model Fig. 4(bottom). For $w_0 < 0$, the model reveals the novel excitability properties associated to a high calcium conductance Fig. 4(top). In the Supplementary Section \ref{SEC: hybrid SH bif}, we further investigate the hybrid transcritical phase portrait though geometrical singular perturbation arguments.

\begin{figure}
\centering
\includegraphics[width=0.8\linewidth]{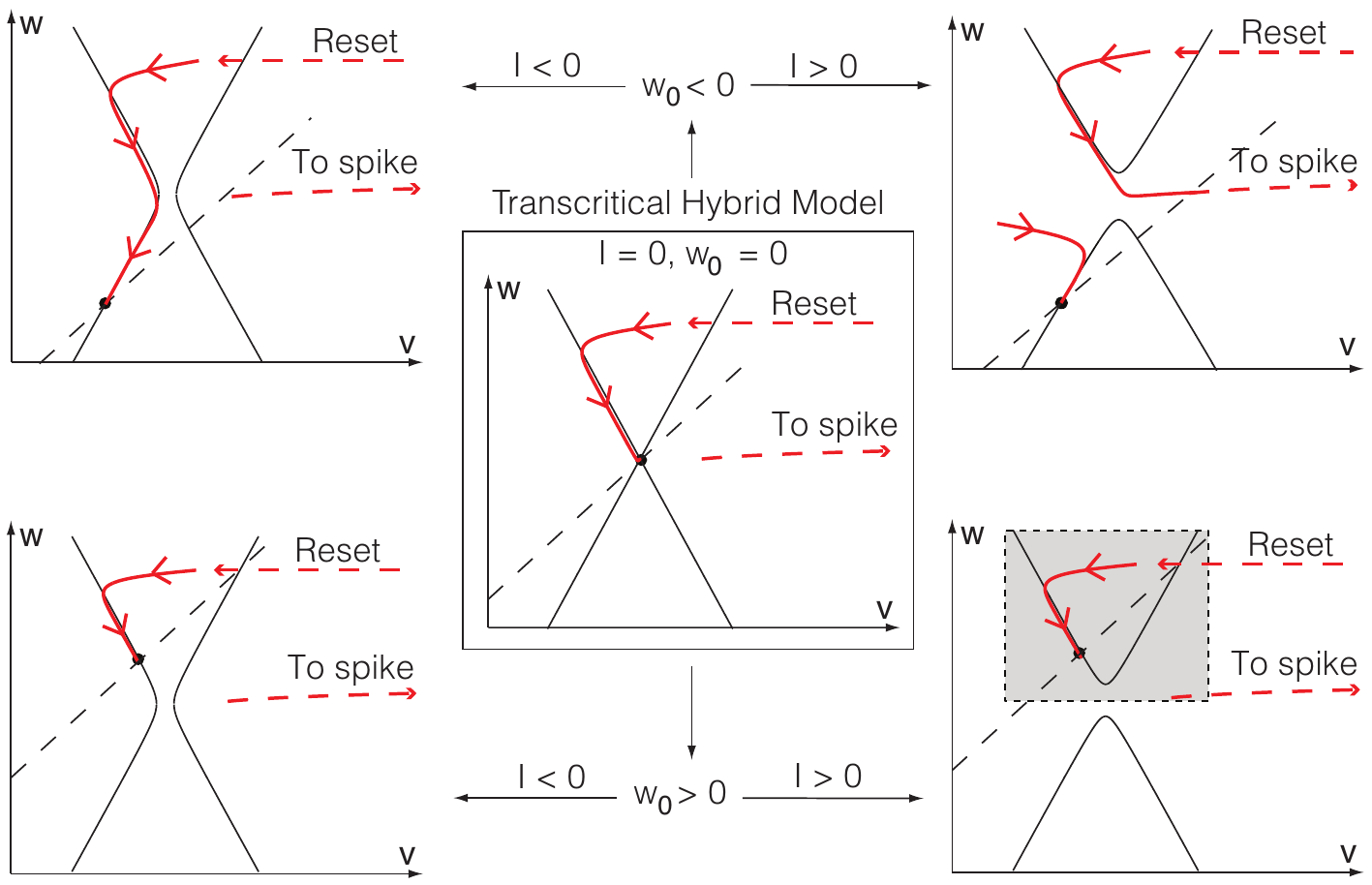}
\caption{{\bf Schematic phase-portraits of the transcritical hybrid model for different values of $I$ and $w_0$.} The $v$- and $n$-nullclines are drawn as full and dashed lines, respectively. The trajectories are drawn as red oriented lines. Many different phase-portraits derive from the transcritical hybrid model, including the one of the fold hybrid model, which only captures the shaded area.}\label{FIG:trans2}
\end{figure}

\section{Hybrid modeling of a thalamocortical relay neuron}
\label{SEC: TC hybrid modeling}

Thalamocortical (TC) relay neurons are the input to sensory cortices. As illustrated in Figure 5, these neurons exhibit two distinct firing patterns: either a continuous regular spiking Fig. 5{\bf a}, or a plateau burst spiking Fig. 5{\bf b}. The switch between the two modes is regulated by prominent T-type calcium currents that are deinactivated by hyperpolarization, thereby modulating the resting membrane potential \cite{main}{McCormick92a}. Thalamocortical neurons are widely studied in the literature, and their two spiking modes make them a prototypical example to illustrate the relevance of the reduced model. We emphasize that our objective is not a fine tuned quantitative modeling of the TC neuron firing pattern. Rather, we attempt to provide a qualitative picture of how the proposed simple hybrid dynamics permits to explain the behavior of TC neurons and, in particular, the role of calcium currents.

Figure 5 compares the experimental step response of a TC neuron {\it in vitro} and the simulated step response of the transcritical hybrid model (\ref{EQ: quantitative IZHw2 model}), both in the low and high calcium conductance modes. As discussed in Section \ref{SEC: new reduced hybrid model}, the small calcium conductance mode is obtained by choosing a positive $w_0$, whereas the large calcium conductance mode is obtained by choosing a negative $w_0$, all the other parameters being identical in the two modes. The hybrid model reproduces the experimental observation: in the low-calcium mode, it responds with a slow regular train of action potentials; in the high-calcium mode, it responds with a long spike latency, plateau oscillations, and an ADP. Phase portrait analysis permits to explain the observed behavior (Supplementary Section S\ref{SEC: hybrid TC phase-plane}). Among other, it shows that the switch from the upper branch of the $v$-nullcline (low calcium conductance mode) to its lower branch (high calcium conductance mode) is sufficient to reproduce and explain the strongly different step-responses of the TC neuron.

\begin{figure}
\centering
\includegraphics[width=0.6\linewidth]{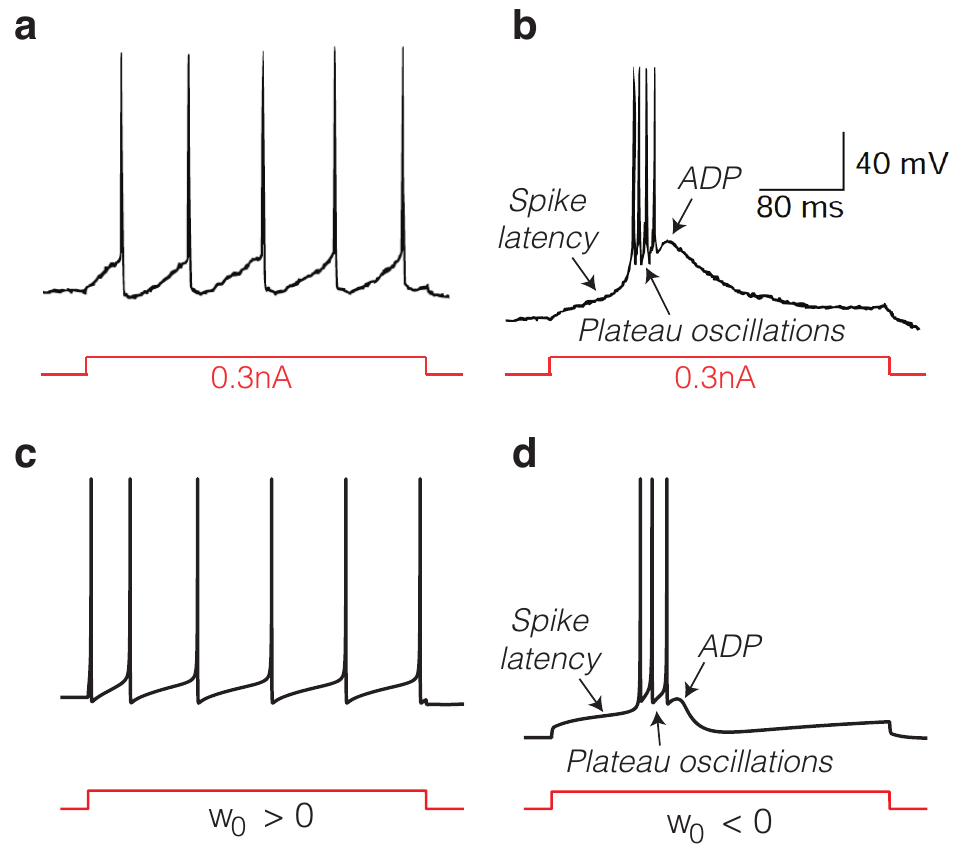}
\caption{{\bf Comparison of the experimental step response of a TC neuron in vitro\cite{main}{Sherman01a} (top) and the step response of the proposed transcritical hybrid model (bottom) in low (left) and high calcium conductance modes (right).} ({\bf a} and {\bf b}): Membrane potential variations of the recorded TC neuron over time in both conditions. ({\bf c} and {\bf d}) Membrane potential variations of the modeled TC neuron over time in both conditions. A variation of $w_0$, which is an image of the calcium conductance, is sufficient to generate the switch of firing pattern physiologically observed in TC cells.}\label{FIG: TC vitro hybrid comp}
\end{figure}

In order to verify the physiological consistence of the transcritical model, we further compare its behavior with the simulated step response of a quantitative 200-compartments model of a TC relay cell\footnote{Simulations were run in the \emph{Neuron} environment, based on the configuration files freely available at {\tt http://cns.iaf.cnrs-gif.fr/alain\_demos.html}.} \cite{main}{destexhe1998dendritic} in the large conductance mode (Fig.6). For the quantitative model, we plot the trajectory projection on the $V-d$ plane, where $V$ and $d$ denotes the somatic membrane potential and the activation gating variable of the somatic T-type calcium current, respectively.

\begin{figure}
\centering
\includegraphics[width=0.8\linewidth]{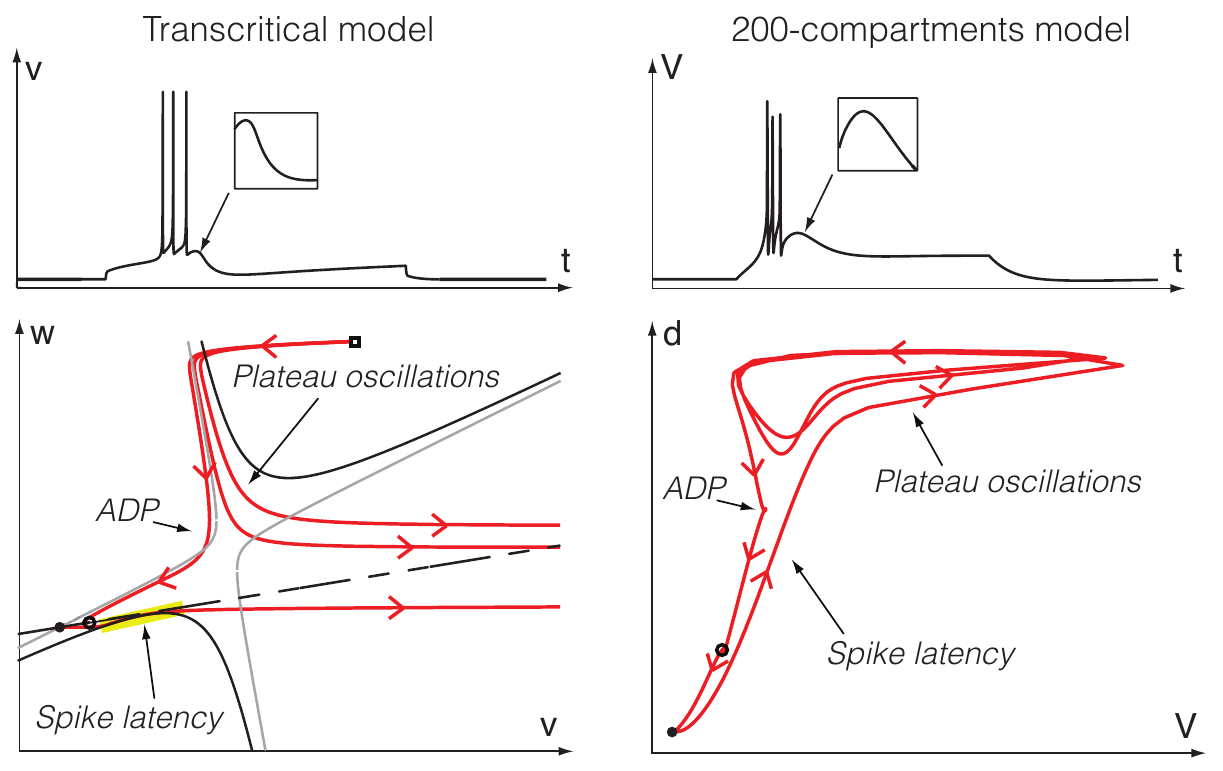}
\caption{{\bf Comparison of the step response and phase-portrait in our proposed hybrid model and in a 200-compartment model of TC neurons \cite{main}{destexhe1998dendritic}, both in large calcium conductance mode.} In the phase-portrait of the transcritical hybrid model, $v$- and $w$-nullclines are drawn as full and dashed lines, respectively, and trajectories are drawn as red oriented lines. The black full line represents the $v$-nullcline at the onset of the stimulation. The gray full line represents the $v$-nullcline at the end of the burst. The phase-portrait of the compartmental model depicts the trajectory projection on the $V-d$ plane, where $V$ and $d$ denotes the somatic membrane potential and the activation gating variable of the somatic T-type calcium current, respectively. In both cases, the ADP is generated during a decrease of the activation variable, and plateau oscillations are exhibited far from the resting state.}\label{FIG: TC hybrid PP}
\end{figure}

As shown on the figure, there is a striking similarity between the projection of the high dimension trajectory and the phase portrait of the second-order transcritical model. In both cases, the ADP is generated during a decrease of the activation variable, and plateau oscillations are exhibited far from the resting state. Moreover, the spike latency is a robust property of the transcritical model because the trajectory must visit the neighborhood of both the nullclines $\dot V = 0$ and $\dot w = 0$ before converging to the spiking limit cycle. It should be stressed that there are no comparable ways to reproduce this behavior in a fold hybrid model. As shown by Izhikevich \cite{main}{IZHIKEVICH2007}, reproducing this behavior with the standard reduced HH model necessitates a non physiological alteration of the reset rule (see also Supplementary Sections \ref{EQ: fold hybrid model} and \ref{SEC: robust ADPs}). This underlines the importance of the revisited model to capture the richness of neuron excitability.

\section{Conclusion}
The inclusion of calcium channels in Hodgkin-Huxley model has a dramatic impact on its mathematical reduction: the firing mechanisms are governed by the local normal form of a transcritical rather than fold bifurcation. Interestingly, it is not the phase portrait of the reduced HH model that is affected by calcium, but only the subregion of the plane where it is physiologically relevant. As a consequence, the classical FitzHugh Nagumo phase portrait is a particular (because localized) view of the more complete picture studied in the present paper.

Although this enlarged phase portrait is the source of rich and diverse forms of excitability, its essence is captured in a simple and physiologically grounded hybrid model. The illustration of its modeling power on the thalamocortical neuron excitability shows the impact of revisiting the classical view. This illustration is just the top of the iceberg because the same principle will apply to many important families of neurons that are thoroughly studied and that have so far largely resisted reduced modeling. The proposed model will impact the understanding of excitability of e.g. dopaminergic, serotonergic, and subthalamic nucleus neurons, whose various firing patterns have a direct and critical impact in physiology and diseases, such as Parkinson's disease and depression. Because of its simplicity and computational efficiency, it is also an ideal candidate for physiologically realistic studies of high-dimensional neuronal networks.

\section{Methods}
\subsection{Equation and parameters of the complete model}
The augmented HH model reads
\begin{eqnarray}\label{EQ: fullHHCA}
C \dot V &=& \overbrace{- \bar g_Kn^4(V-V_K) - \bar g_{Na}m^3h(V-V_{Na})- g_l(V-V_l)+I_{app}}^{\text{Hodgkin-Huxley dynamics}} \overbrace{+ I_{Ca} + I_{pump}}^{\text{Calcium currents}}\nonumber\\
\dot n&=&\alpha_n(V)(1-n)-\beta_n(V)n \nonumber\\
\dot m&=&\alpha_m(V)(1-m)-\beta_m(V)m \nonumber\\
\dot h&=&\alpha_h(V)(1-h)-\beta_h(V)h, \nonumber
\end{eqnarray}

For the HH dynamics, we use the parameters of the original paper \cite{main}{HODHUX}. As all other HH currents, the additional calcium current obeys Ohm's law
$$I_{Ca}=-\bar g_{Ca}d^a(V-V_{Ca})$$

where $\bar g_{Ca}$ is the maximum calcium conductance, $V_{Ca}$ is the calcium Nernst potential, and $d$ is the calcium activation gating variable. Such kinetics fit the majority of calcium channel subtypes \cite{main}{catterall2005international}. Exploiting the similarities between the potassium and calcium gating kinetics, we further assume that the behavior of the calcium activation gating variable $d$ is well approximated as a static function of the potassium activation gating variable $n$. The simple choice $d:=n$, and $a:=3$ suffices for our paper. The exponent accounts for the delayed activation of calcium currents, described often empirically by second to sixth powers (see \cite{main}{hille1991ionic} page 112, and references therein). These parameters do not reflect any precise physiological calcium current. We choose them as a prototypical example. The functions $\alpha_x,\beta_x$, $x=n,m,h$, can be found in the paper \cite{main}{HODHUX}. The value for the potassium Nernst potential $V_K=-12$ is the same as in \cite{main}{HODHUX}, while the sodium Nersnt and the leak Nernst potential are rounded to $V_{Na}=120$ and $V_l=10.6$, respectively. The values of the sodium $\bar g_{Na}=120$, potassium $\bar g_K=36$, and leak $g_l$ (maximum) conductances are the same as in \cite{main}{HODHUX}. The calcium Nernst potential is given by $V_{Ca}=150$. The numerical simulations of Figure 1(right) are obtained by picking $g_{Ca}=2.7$ and $I_{pump}=-17$. Parameter unchanged for Equations (\ref{EQ: redHHCA}) and Figure \ref{FIG: global phase portrait}.

\subsection{Planar reduction and phase portrait analysis}
We follow the standard reduction of the original HH model to a two dimensional system by: i) assuming an instantaneous sodium activation, $m\equiv m_\infty(V)$, where $m_\infty(V)=\alpha_m(V)/(\alpha_m(V)+\beta_m(V))$; ii) exploiting the approximate linear relation, originally proposed in \cite{main}{fitzhugh61}, $h+sn=c$, with $s,c\simeq 1$. Applying the same reduction to (\ref{EQ: fullHHCA}) with parameters as above, we obtain the planar system
\begin{IEEEeqnarray}{rCl}\label{EQ: redHHCA}
C\dot V &=& - \bar g_Kn^4(V-V_K) - \bar g_{Na}m_\infty(V)^3(0.89-1.1n)(V-V_{Na})-g_l(V-V_l) + I_{app} \nonumber\\
	&&- \bar g_{Ca}n^3(V-V_{Ca}) + I_{pump}\IEEEyessubnumber\\
\dot n&=&\alpha_n(V)(1-n)-\beta_n(V)n \IEEEyessubnumber
\end{IEEEeqnarray}

$V$-``nullcline'' refers to the set $\{(V,n)\in\mathbb R^2:\ \dot V=0\}$ and similarly for other variables.

\subsection{Hybrid modeling of TC neurons}
For modeling convenience, we add a fitting parameter $b\in\mathbb R$ in the sub-threshold (continuous) voltage dynamics: 
$$\dot v=v^2+bvw-w^2+I$$

The extra parameter  does not affect the nature of the bifurcation, but tunes the slope of the $v$-nullcline branches. In order to account for the effect of intracellular calcium variations and for the associated activation of calcium pump currents, we also add a slow adaptation variable $z$ as follows:
\begin{IEEEeqnarray}{rCl}\label{EQ: quantitative IZHw2 model}
\dot v=v^2+bvw-w^2+I-z&\quad\quad\quad \mbox{if }\ v\geq v_{th},\ \mbox{then}\IEEEyessubnumber\\
\dot w=\epsilon(av-w+w_0)&\quad\quad\quad v\leftarrow c,\ w\leftarrow d,\IEEEyessubnumber\\
\dot z=-\epsilon_zz,&\quad\quad\quad z\leftarrow z+d_z\IEEEyessubnumber.
\end{IEEEeqnarray}

with $a=0.1,\,b=-3,\,c=15,\,d=15,\,\epsilon=1,\,v_{th}=100,\,\epsilon_z=0.1,\,d_z=40$. The low-calcium mode correspond to $w_0=3.2$, the high-calcium mode to $w_0=-4$. The injected current $I=-5$ when the stimulation is off and $I=85$ when the stimulation is on. We stress two properties that are not captured by the reduced model (\ref{EQ: quantitative IZHw2 model}): Firstly, the intracellular calcium ({\it i.e.} $z$) dynamics are generally coupled with the membrane voltage also in the subthreshold phase. Secondly, depending on the type of the modeled calcium current, the calcium conductance, here reflected by $w_0$ (Section \ref{SEC: new reduced hybrid model}), might slowly change in response to voltage and intracellular calcium variations. Both properties are neglected in (\ref{EQ: quantitative IZHw2 model}).

\subsection{Numerical simulations}
Numerical simulations were run with MATLAB \\
({\tt http://www.mathworks.com}), apart from the 200-compartment model in Figures 5 and6 which was simulated with the NEURON software environment\\ 
({\tt http://www.neuron.yale.edu}). The numerical bifurcation analysis of Supplementary Section S2 has been obtained with XPP environment\\
({\tt http://www.math.pitt.edu/$\thicksim$bard/xpp/xpp.html})



\section*{Acknowledgments}
The research leading to these results has received funding from the European Union Seventh Framework Programme [FP7/2007-2013]  under grant agreement n$^\circ$257462 HYCON2 Network of excellence, by grant 9.4560.03 from the F.R.S.-FNRS (VS), and by two grants from the Belgian Science Policy (IAP6/31 (VS) and IAP6/4 (RS)).

\bibliographystyle{main}{plain}
\bibliography{main}{../../../refs}{References}

\clearpage

\begin{center}
{\huge Supplementary material}
\end{center}


\newbibliography{supp}

 
\renewcommand{\thesection}{S.\arabic{section}}
\renewcommand{\thesubsection}{\thesection.\arabic{subsection}}
 
%
\makeatletter 
\def\tagform@#1{\maketag@@@{(S\ignorespaces#1\unskip\@@italiccorr)}}
\makeatother
 
\makeatletter
\makeatletter \renewcommand{\fnum@figure}
{\figurename~S\thefigure}
\makeatother




\setcounter{section}{0}
\setcounter{figure}{0}

\renewcommand{\figurename}{Figure}

\section{Supplementary figures}

\begin{figure}[h!]
\centering
\includegraphics[width=0.9\textwidth]{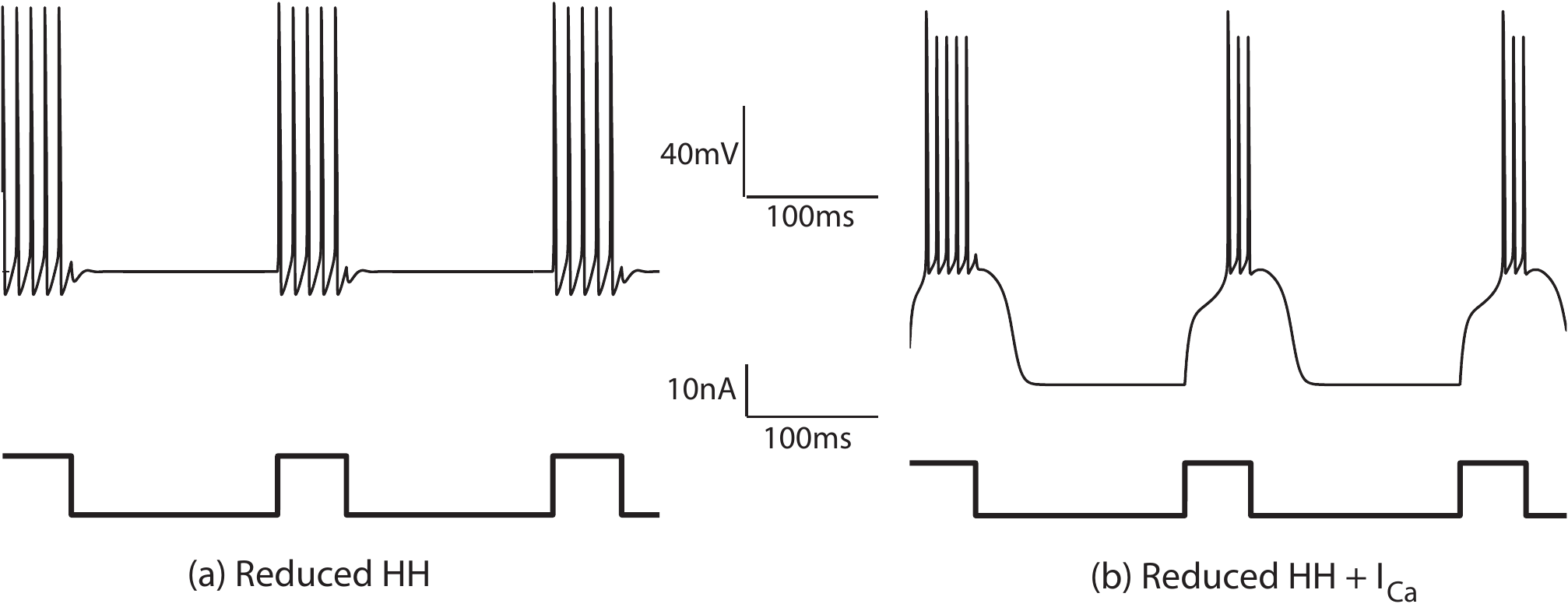}
\caption{{\bf Step responses of the HH model without (left) and with a calcium current (right).} Time-evolution of the applied excitatory current (bottom) and of the corresponding membrane potential (top).}\label{FIG: red comparison}
\end{figure}

\begin{figure}[h!]
\center
\includegraphics[width=0.8\textwidth]{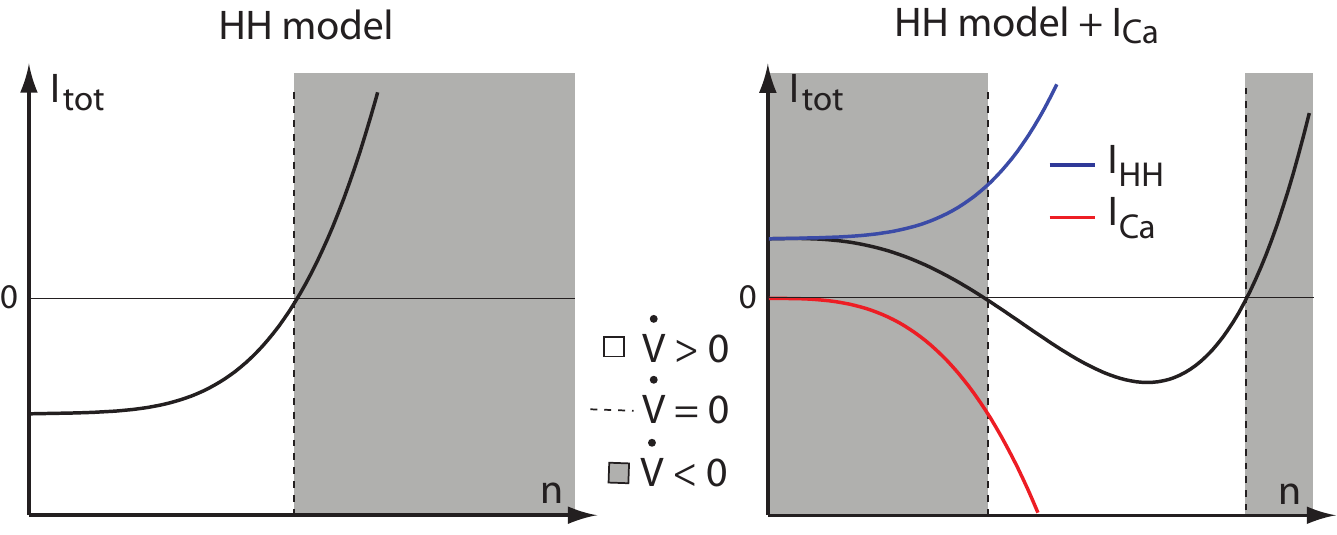}
\caption{{\bf Total ionic currents for $V$ fixed as a function of $n$ without (left) and with calcium channels (right).} Blank portions corresponds to the values of $n$ where $\dot{V} > 0$, shaded portions corresponds to the values of $n$ where $\dot{V} < 0$, and the dashed lines correspond to the values of $n$ where $\dot{V} = 0$. Note that the total ionic currents monotonically increase only in the absence of calcium channels (left).}\label{EQ: V-null phys interpret}
\end{figure}

\begin{figure}[h!]
\center
\includegraphics[width=0.8\textwidth]{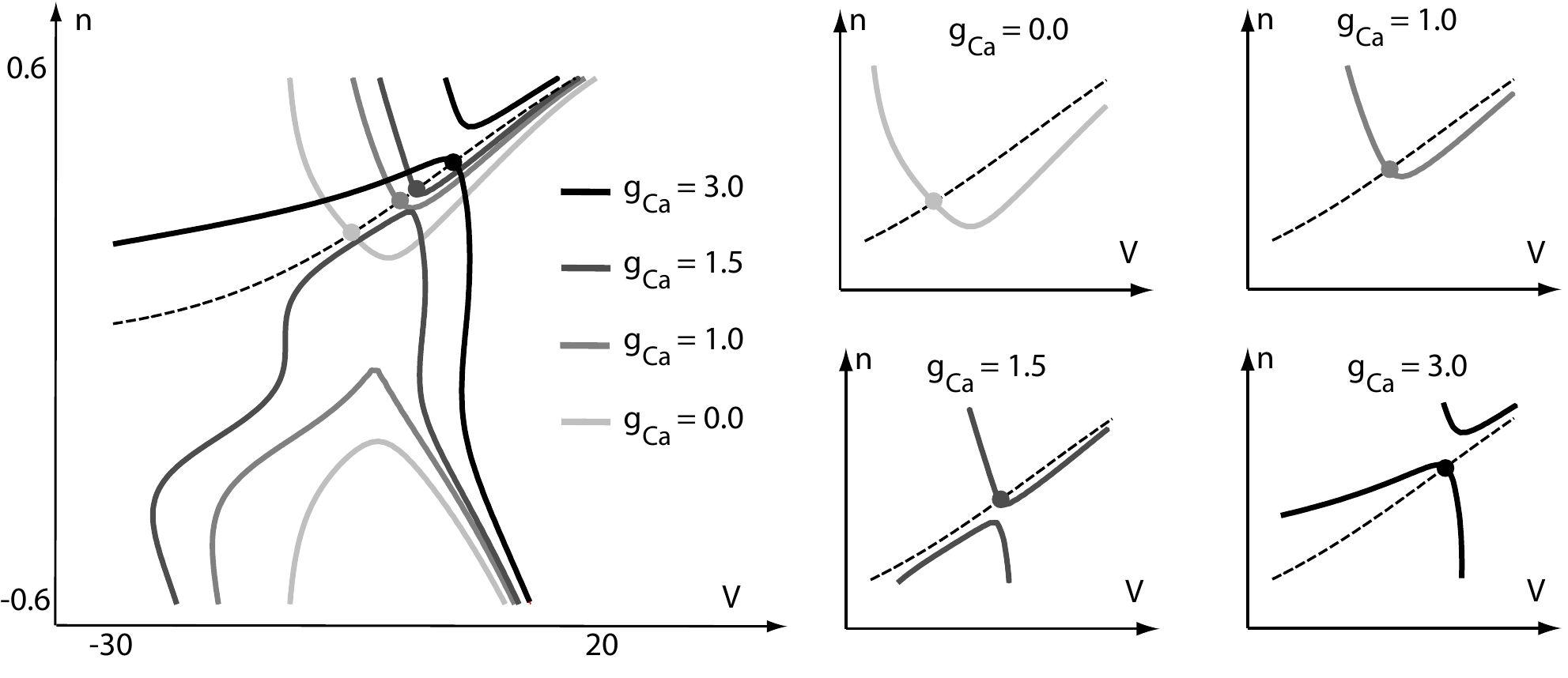}
\caption{{\bf Nullcline intersections in the reduced model (5) with calcium current for different values of the calcium conductance.} The $n$-nullcline is depicted as a dashed line, the $V$-nullcline as a solid line. The associated calcium conductance is expressed via a gray scale, as indicated in the figure legend. The calcium pump currents is given by $I_{Ca,pump}=0.0,\,-0.74,\,-2.78,\,-15.6$ for, respectively, $g_{Ca}=0.0,\,1.0,\,1.5,\,3.0$.}\label{FIG: Na to Ca}
\end{figure}

\clearpage

\section{Bifurcation analysis}
\label{SEC: bif analysis}

A bifurcation diagram of (5) with $I_{app}$ as the bifurcation parameter sheds more light on the transition mechanism between the resting and spiking modes (Figure \ref{FIG:1}). We use XPPAUT \cite{supp}{ERMEXPPbook} for this numerical analysis. In both the $I_{Ca}$-on and off modes, we draw the bifurcation diagram only for small $I_{app}$, corresponding to the transition from resting to limit cycle oscillations (for larger $I_{app}$, the stable limit cycle disappears in a supercritical Andronov-Hopf bifurcation in both cases, which leads to a stable depolarized, {\it i.e.} high-voltage, state).

\begin{figure}[h!]
\centering
\includegraphics[width=0.9\textwidth]{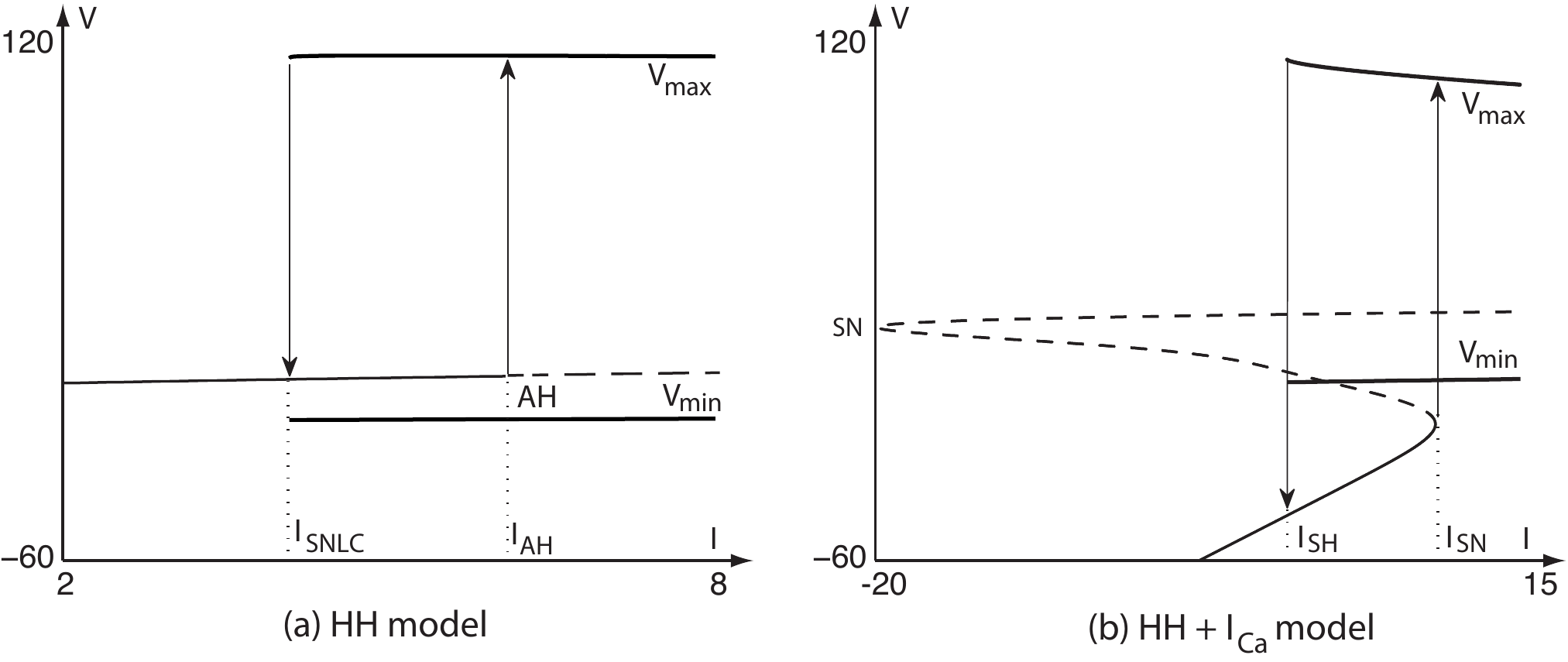}
\caption{{\bf One parameter bifurcation diagram of the reduced Hodgkin-Huxley model (5) with and without calcium current.} Thin solid lines represents stable fixed points, while dashed lines unstable fixed points or saddle points. The thick lines labeled $V_{min}$ and $V_{max}$ represent the minimum and the maximum voltage of stable limit cycles, respectively. (a): Hodgkin-Huxley model. (b): Hodgkin-Huxley model with calcium current. $I_x$, with $x=SNLC,AH,SH,SN$, denotes the value of the input current for which the system undergoes the bifurcation $x$. More details are given in the text.
}
\label{FIG: red bif comparison} 
\end{figure}

Figure S\ref{FIG: red bif comparison} (left) illustrates the bifurcation diagram of the original reduced Hodgkin-Huxley model. For low values of $I_{app}$, the unique fixed point is a stable focus that loses stability in a subcritical Andronov-Hopf bifurcation at $I_{app}=I_{AH}$. Beyond the bifurcation, the trajectory converges to the stable spiking limit cycle. When $I_{app}$ is lowered again below $I_{SNLC}$, the spiking limit cycle disappears in a saddle-node of limit cycles, the unstable one (not drawn) emanating from the subcritical Andronov-Hopf bifurcation, and the trajectory relaxes back to rest.


Figure S\ref{FIG: red bif comparison} (right) illustrates the bifurcation diagram of the reduced Hodgkin-Huxley model in the presence of calcium channels. For $I_{app}<I_{SN}$, a stable node (lower branch), a saddle (central branch), and an unstable focus (upper branch) are present, as in Figure \ref{FIG:1}{\bf b}(top right). The node and the saddle coalesce in a supercritical fold bifurcation at $I_{app}=I_{SN}$, and disappear for $I_{app}>I_{SN}$, letting the trajectory converge toward the stable limit cycle. The spike latency observed in the $I_{Ca}$-on configuration unmasks the ghost of this bifurcation. The stable limit cycle disappears in a saddle homoclinic bifurcation as $I_{app}$ falls below $I_{SH}$, which lets the trajectory relax back to the hyperpolarized state.

The homoclinic bifurcation exhibited by the Hodgkin-Huxley model with calcium channels is a key mathematical difference with respect to the standard HH model. The next section unfolds this bifurcation by exploiting the time-scale separation between the fast voltage $V$ and the slow recovery variable $n$.

\clearpage
\section{Transcritical and saddle-homoclinic bifurcations of in the presence of calcium channels}
\label{SEC: trans and homo bifurcations}

By exploiting the sharp time-scale separation between the membrane voltage and the recovery variable dynamics, in this section we assume that $\dot n=O(\epsilon)\dot V$, where $\epsilon>0$, and study the singularly perturbed limit $\epsilon=0$. We highlight in particular the existence of a transversal $V$-nullcline self-intersection in the passage from Figure \ref{FIG:1}{\bf b}(top right) to Figure \ref{FIG:1}{\bf b}(bottom right) that is associated to a transcritical singularity (see also Fig.~\ref{FIG: trans}{\bf a} center). We then rely on some recent results \cite{supp}{KRSZ01} to qualitatively study the phase portrait of (\ref{EQ: redHHCA}) away from the singular limit ({\it i.e.} $\epsilon>0$) and for different values of the input current.

\subsection{Transcritical bifurcation in the presence of calcium channels}
\label{SEC: trans bif}
%

\begin{figure}
\centering
\subfigure[][$I=2.4$]{
\includegraphics[width=0.45\textwidth]{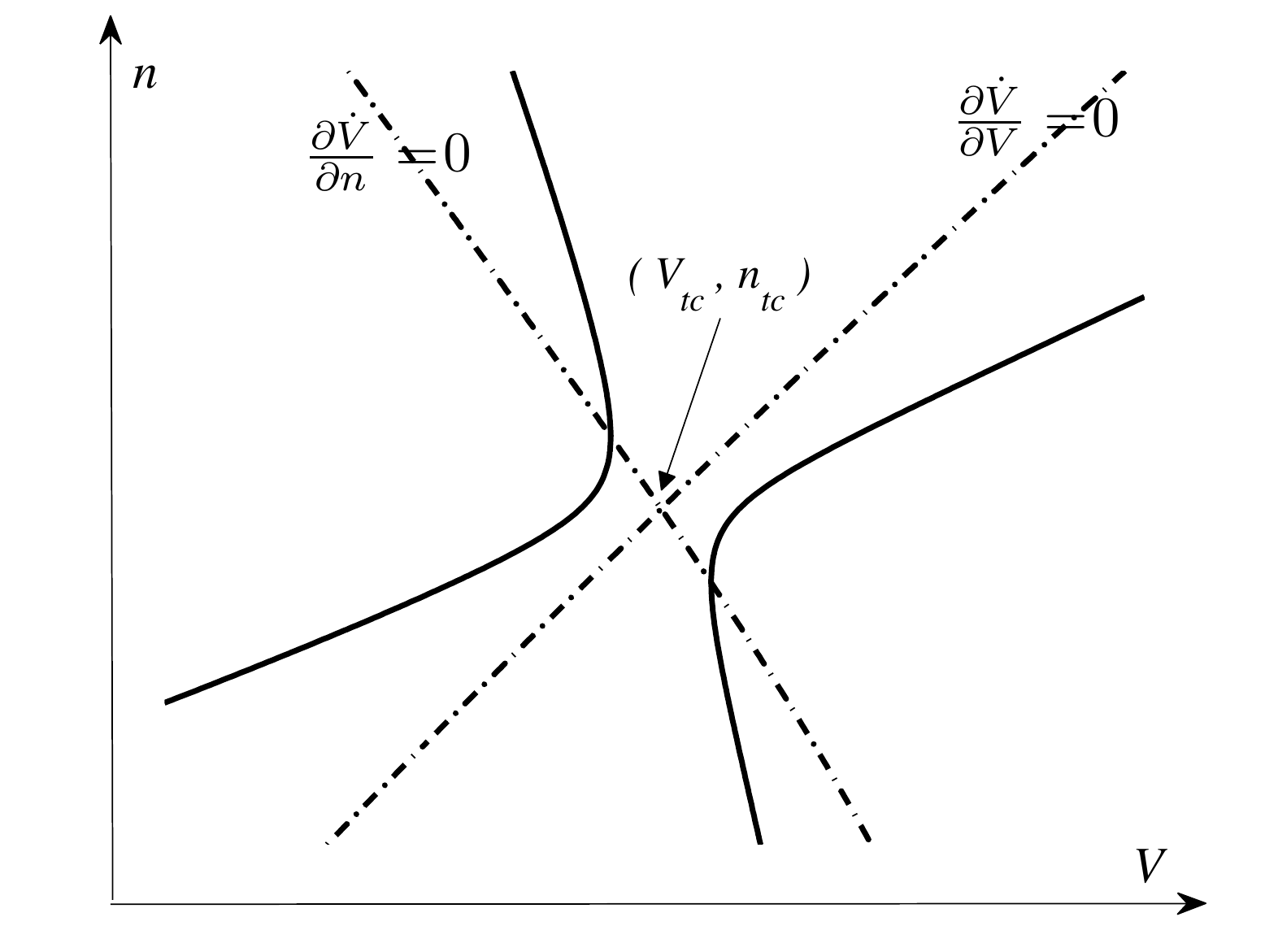}
}
\subfigure[][$I=2.46$]{
\includegraphics[width=0.45\textwidth]{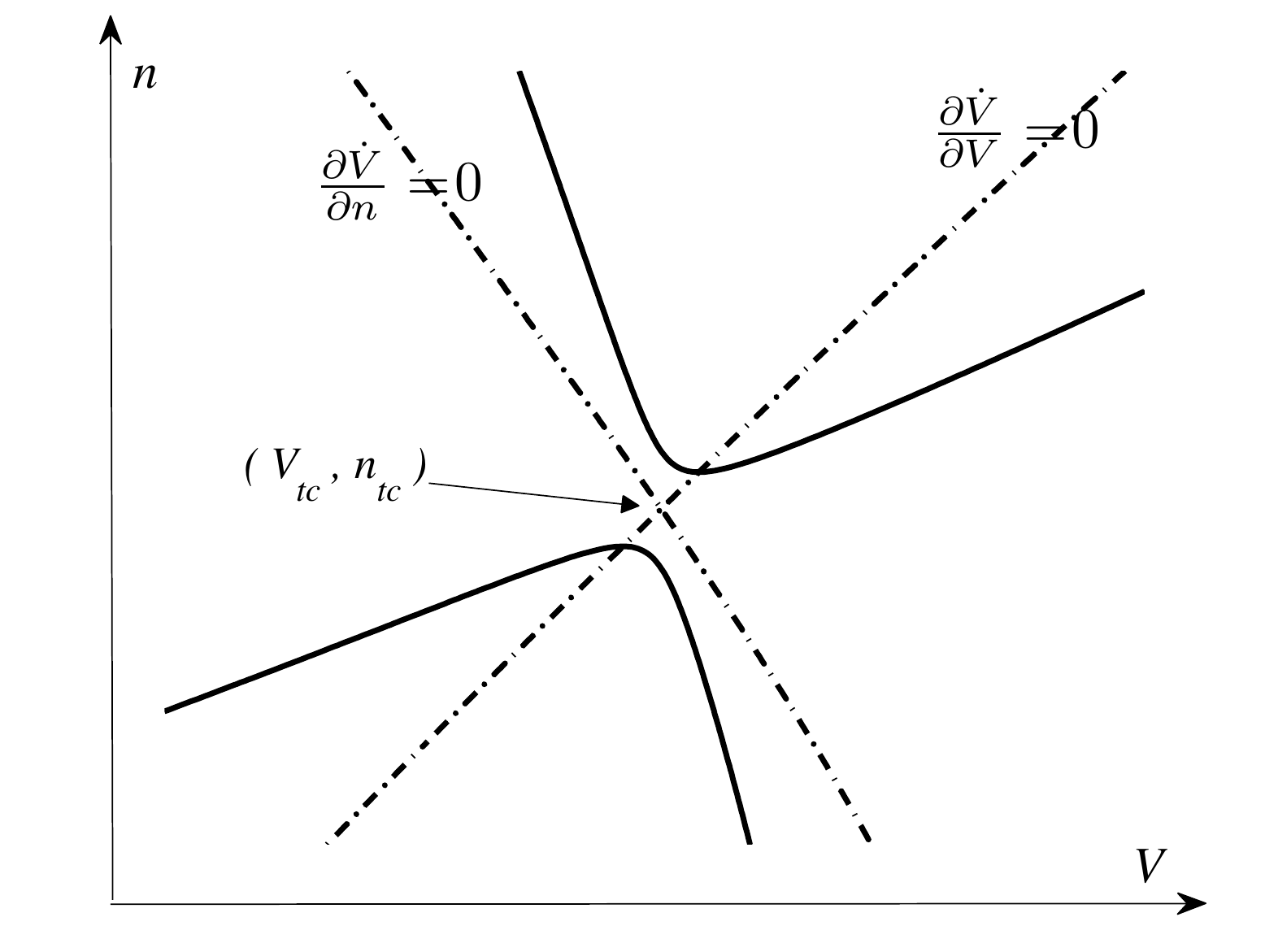}
}
\caption{{\bf Phase portrait of the model for different values of the input current in the presence of calcium channels.} The $V$-nullcline is drawn as a solid line. The dashed-lines depict the locus $\frac{\partial f_V}{\partial V}=0$ and $\frac{\partial f_V}{\partial n}=0$, respectively.}
\label{FIG: zoomed phase portrait} 
\end{figure}

The presence of the $V$-nullcline self intersection is confirmed by a closer inspection of the phase portrait near the nullcline break-up, as illustrate in Figure S\ref{FIG: zoomed phase portrait}. Let us discuss this transition in more details. The existence of the self intersection follows directly by $\frac{\partial\dot V}{\partial I_{app}}=1$ and the implicit function theorem \cite{supp}{Lee_v2006}.
We claim that it is also transversal. Let $I_{tc}\in (2.4,2.46)$ be the input current for which the $V$-nullcline has the self-intersection. Let $f_V(V,n,I_{app})$ be the voltage dynamics of the reduced model (\ref{EQ: redHHCA}) in the presence of calcium current, that is
\begin{equation}\label{EQ: vdot trans cond}
f_V(V,n,I_{app}):=\dot V\big |_{g_{Ca}=2.7,\ I_{pump}=-17}\quad,
\end{equation}
where $\dot V$ is defined in (\ref{EQ: redHHCA}) and all the other parameters are unchanged. Then, we want to show that there exists $(V_{tc},n_{tc})\in\mathbb R^2$, such that
\begin{subequations}\label{EQ: trans cond}
\begin{eqnarray}
f_V(V_{tc},n_{tc},I_{tc})&=&0\label{EQ: trans cond a}\\
\frac{\partial f_V}{\partial V}(V_{tc},n_{tc},I_{tc})&=&0\label{EQ: trans cond b}\\
\frac{\partial f_V}{\partial n}(V_{tc},n_{tc},I_{tc})&=&0\label{EQ: trans cond c}\\
\left|
\begin{array}{cc}
\displaystyle{\frac{\partial^2 f_V}{\partial V^2}}(V_{tc},n_{tc},I_{tc}) & \displaystyle{\frac{\partial^2 f_V}{\partial V\partial n}}(V_{tc},n_{tc},I_{tc})\vspace{2mm}\\
\displaystyle{\frac{\partial^2 f_V}{\partial V\partial n}}(V_{tc},n_{tc},I_{tc}) & \displaystyle{\frac{\partial^2 f_V}{\partial n^2}}(V_{tc},n_{tc},I_{tc})
\end{array}
\right|&<&0\label{EQ: trans cond d}\\
\frac{\partial^2 f_V}{\partial V^2}(V_{tc},n_{tc},I_{tc})&\neq&0,\label{EQ: trans cond e}
\end{eqnarray}
\end{subequations}
describing a non-degenerate (transversal) self-intersection of the $V$-nullcline at $(V_{tc},n_{tc})$ (see e.g.  \cite{supp}{KRSZ01} and Section 5.5.2, Th. 5.7 in \cite{supp}{SEYDEL94}). The self-intersection satisfies the additional constraint that both intersecting branches are not parallel to the $V$ axis, as implied by (S\ref{EQ: trans cond e}). Condition (S\ref{EQ: trans cond a}) is the $V$-nullcline equation. Noticing that $\frac{\partial f_V}{\partial n}$ and $\frac{\partial f_V}{\partial V}$ do not depend on $I_{app}$, conditions (S\ref{EQ: trans cond b}) and (S\ref{EQ: trans cond c}) follows by the fact that, as $I_{app}<I_{tc}$ varies, as in Figure S\ref{FIG: zoomed phase portrait}(a), the right and left extrema of, respectively, the left and right branches of the $V$-nullcline lie, by definition, on the line $\frac{\partial f_V}{\partial n}=0$, and, similarly, as $I_{app}>I_{tc}$ varies, as in Figure S\ref{FIG: zoomed phase portrait}(b), the minimum and the maximum of, respectively, the upper and lower branches of the $V$-nullcline lie, by definition, on the line $\frac{\partial f_V}{\partial V}=0$. Since at the intersection the four extrema coincide, conditions (S\ref{EQ: trans cond b}) and (S\ref{EQ: trans cond c}) follow. We stress that (S\ref{EQ: trans cond b}) and (S\ref{EQ: trans cond c}) define the point $(V_{tc},n_{tc})$ as the intersection of two lines (cf. Figure S\ref{FIG: zoomed phase portrait}) in a unique way. Conditions (S\ref{EQ: trans cond d}) and (S\ref{EQ: trans cond e}) are generic and can be easily verified numerically.


In the singular limit $\dot n=0$, the self intersection described by conditions (S\ref{EQ: trans cond}) corresponds to a transcritical bifurcation (see e.g. Section Section 3.2 in \cite{supp}{strogatz00}) of the voltage dynamics, with $n$ as the bifurcation parameter. That is, as sketched in Figure \ref{FIG: trans}{\bf a}(center), for $I_{app} =I_{tc}$, the two intersecting branches exchange their stability at $n=n_{tc}$. More precisely, as $n$ varies above and below the intersection point, {\it i.e.} $n\neq n_{tc}$, the voltage dynamics has two fixed points, the left one is stable, whereas the right is unstable. The stability of the fixed points follows from the fact that, on the left of the self-intersection, $\frac{\partial f_V}{\partial V}<0$, while $\frac{\partial f_V}{\partial V}>0$  on the right (cf. Figure S\ref{FIG: zoomed phase portrait}). The two fixed points collide in a transcritical singularity at the self-intersection.


\subsubsection{A normal form lemma}

In this technical section we compute a normal form of (\ref{EQ: redHHCA}) associated to the self-intersection described by conditions (S\ref{EQ: trans cond}).

Give $\epsilon>0$, let
\begin{equation}\label{EQ: eps-scaled ndot}
\epsilon g(V,n):=\dot n,
\end{equation}
where $\dot n$ is defined in (\ref{EQ: redHHCA}), and
\begin{equation}\label{EQ: ndot at trans}
g_0:=g(V_{tc},n_{tc}),
\end{equation}
where $V_{tc}$ and $n_{tc}$ satisfy the defining conditions (S\ref{EQ: trans cond b})-(S\ref{EQ: trans cond b}). The form (S\ref{EQ: eps-scaled ndot}) is just a rescaling (through $\epsilon$) of the recovery variable dynamics that highlights the timescale separation between $V$ and $n$. The following lemma is an application of Lemma 2.1 in \cite{supp}{KRSZ01} to the reduced HH model (\ref{EQ: redHHCA}) with calcium current.

\begin{lemma}\label{LEM: trans normal form}
Suppose that $g_0<0$. Then there exists an affine change of coordinates and a rescaling of $\epsilon$ and $I$ that transforms the dynamical system
\begin{subequations}\label{EQ: trans dyn syst}
\begin{eqnarray}
\dot V&=&f_V(V,n,I_{app})\IEEEyessubnumber\\
\dot n&=&\epsilon g(V,n)\IEEEyessubnumber,
\end{eqnarray}
\end{subequations}
where $f_V$ and $g$ are respectively defined in (S\ref{EQ: vdot trans cond}) and (S\ref{EQ: eps-scaled ndot}), into the dynamical system
\begin{subequations}\label{EQ: trans normal form}
\begin{eqnarray}
\dot v&=&v^2-w^2+I+h_1(v,w,\epsilon)\IEEEyessubnumber\\
\dot w&=&\epsilon(-1+h_2(v,w,\epsilon))\IEEEyessubnumber,
\end{eqnarray}
\end{subequations}
where $h_1(v,w,\epsilon)=O(v^3,v^2w,vw^2,w^3,\epsilon v,\epsilon w, \epsilon^2)$ and $h_2(v,w,\epsilon)=O(v,w,\epsilon)$.
\end{lemma}
\begin{proof}
Let $\alpha:=\frac{1}{2}\frac{\partial^2 f_V}{\partial V^2}(V_{tc},n_{tc},I_{tc})$, $\beta:=\frac{1}{2}\frac{\partial^2 f_V}{\partial V\partial n}(V_{tc},n_{tc},I_{tc})$, $\gamma:=\frac{1}{2}\frac{\partial^2 f_V}{\partial n^2}(V_{tc},n_{tc},I_{tc})$, and
$\lambda:=\frac{-\beta}{\sqrt{\beta^2-\gamma\alpha}}$.
From (\ref{EQ: trans cond}) and \cite[Lemma 2.1]{supp}{KRSZ01}, it follows that, for $I=I_{tc}$, the affine change of variable $ v=\alpha (V-V_{tc})+\beta (n-n_{tc})$, $w=\sqrt{\beta^2-\gamma\alpha}(n-n_{tc})$, transforms (\ref{EQ: trans dyn syst}), after a suitable rescaling of $\epsilon$, into the equation
\begin{IEEEeqnarray*}{rCl}
\dot v&=&v^2-w^2+\lambda\epsilon+h_1(v,w,\epsilon)\\
\dot w&=&\epsilon(-1+h_2(v,w,\epsilon)).
\end{IEEEeqnarray*}
Noticing that $f_V$ is affine in the input current, the extra term $\alpha(I-I_{tc})$ must be added to $\dot v$ in the case $I\neq I_{tc}$. The result follows by defining the rescaled input current $\tilde I:=\lambda\epsilon+\alpha(I-I_{tc})$.
\end{proof}\\

Note that for the dynamical system (S\ref{EQ: trans normal form}) with $I=0$, the self intersection of the $V$-nullcline discussed in Section \ref{SEC: trans bif} becomes evident, with the two intersecting branches given by $w=\pm v$.

\subsection{Singularly perturbed saddle-homoclinic bifurcation in the presence of calcium channels}
\label{SEC: SP homoclinic bifurcation}

\begin{figure}
\centering
\subfigure[][$I=2.4<I_{SH}$]{
\includegraphics[width=0.45\textwidth]{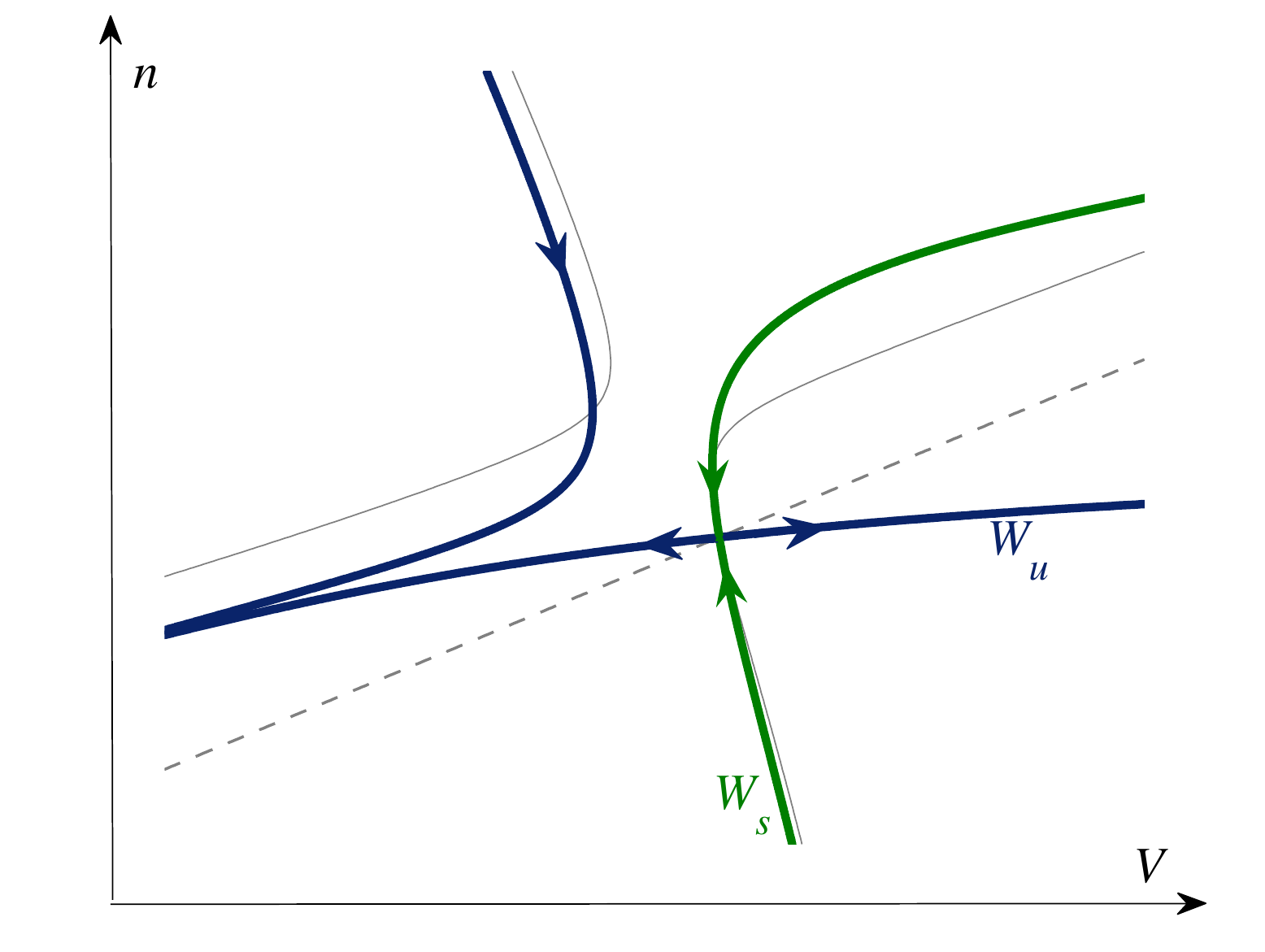}
}
\subfigure[][$I=2.46<I_{SH}$]{
\includegraphics[width=0.45\textwidth]{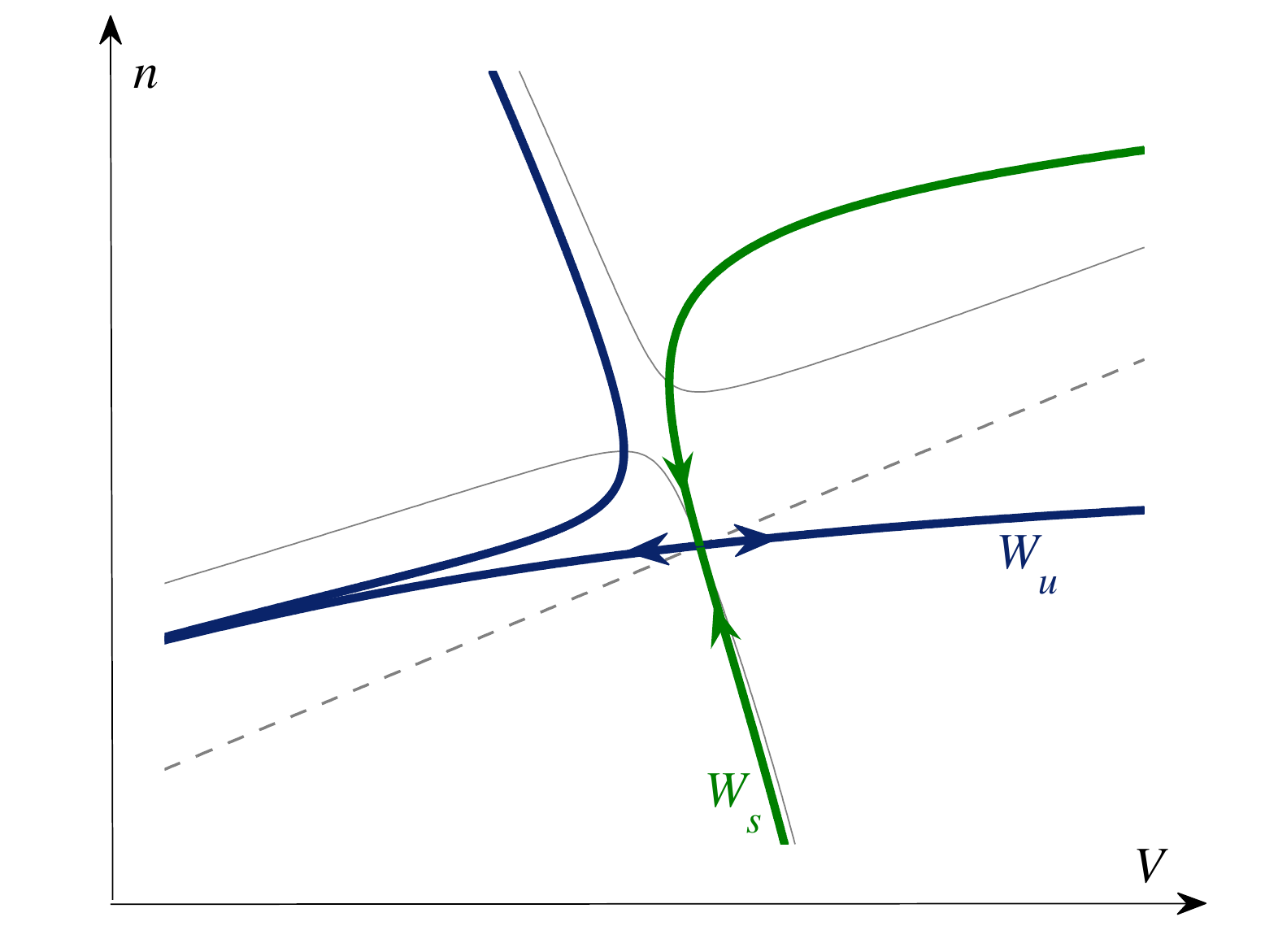}
}
\subfigure[][$I=I_{SH}\sim2.497$]{
\includegraphics[width=0.45\textwidth]{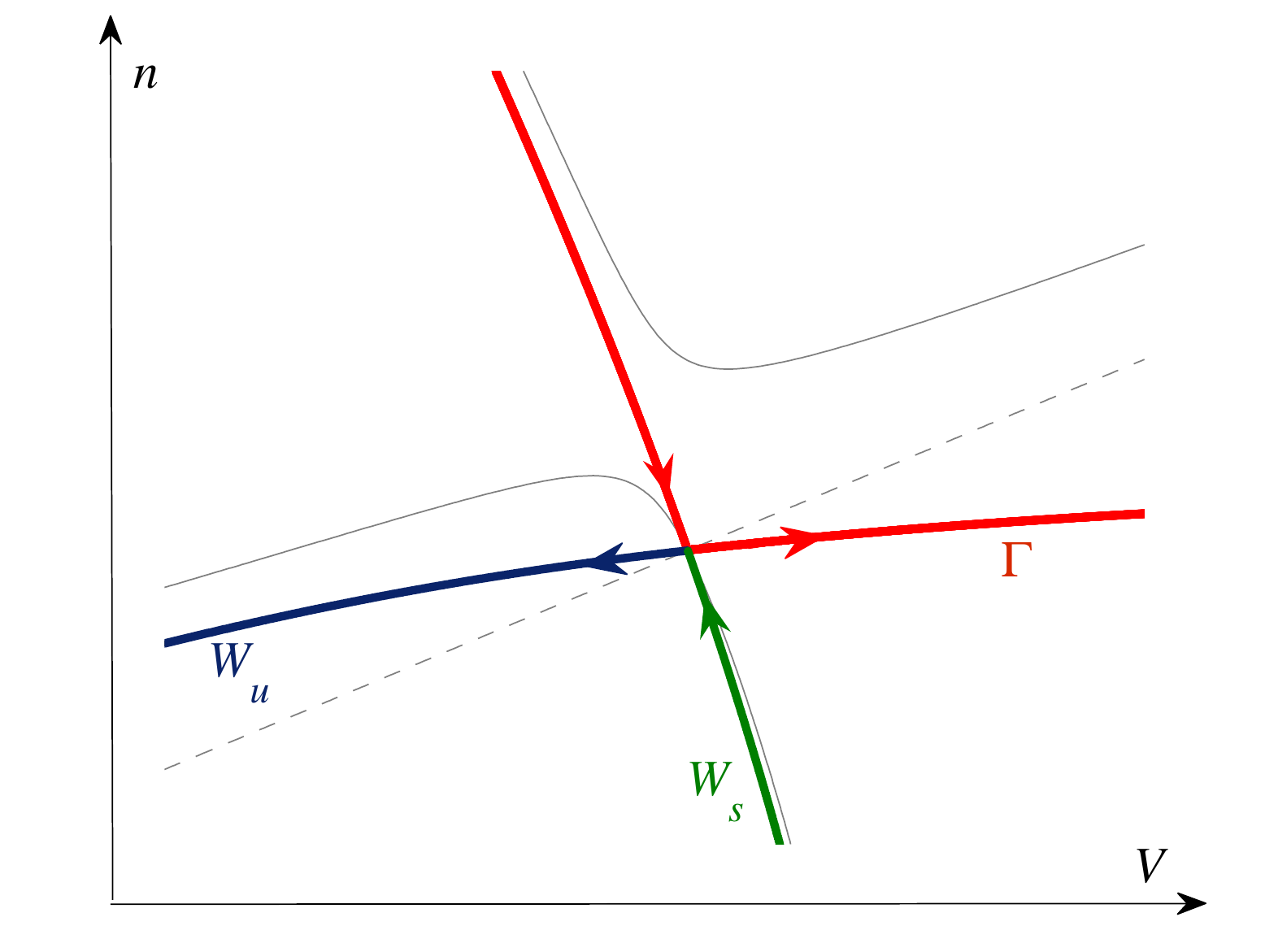}
}
\subfigure[][$I=2.52>I_{SH}$]{
\includegraphics[width=0.45\textwidth]{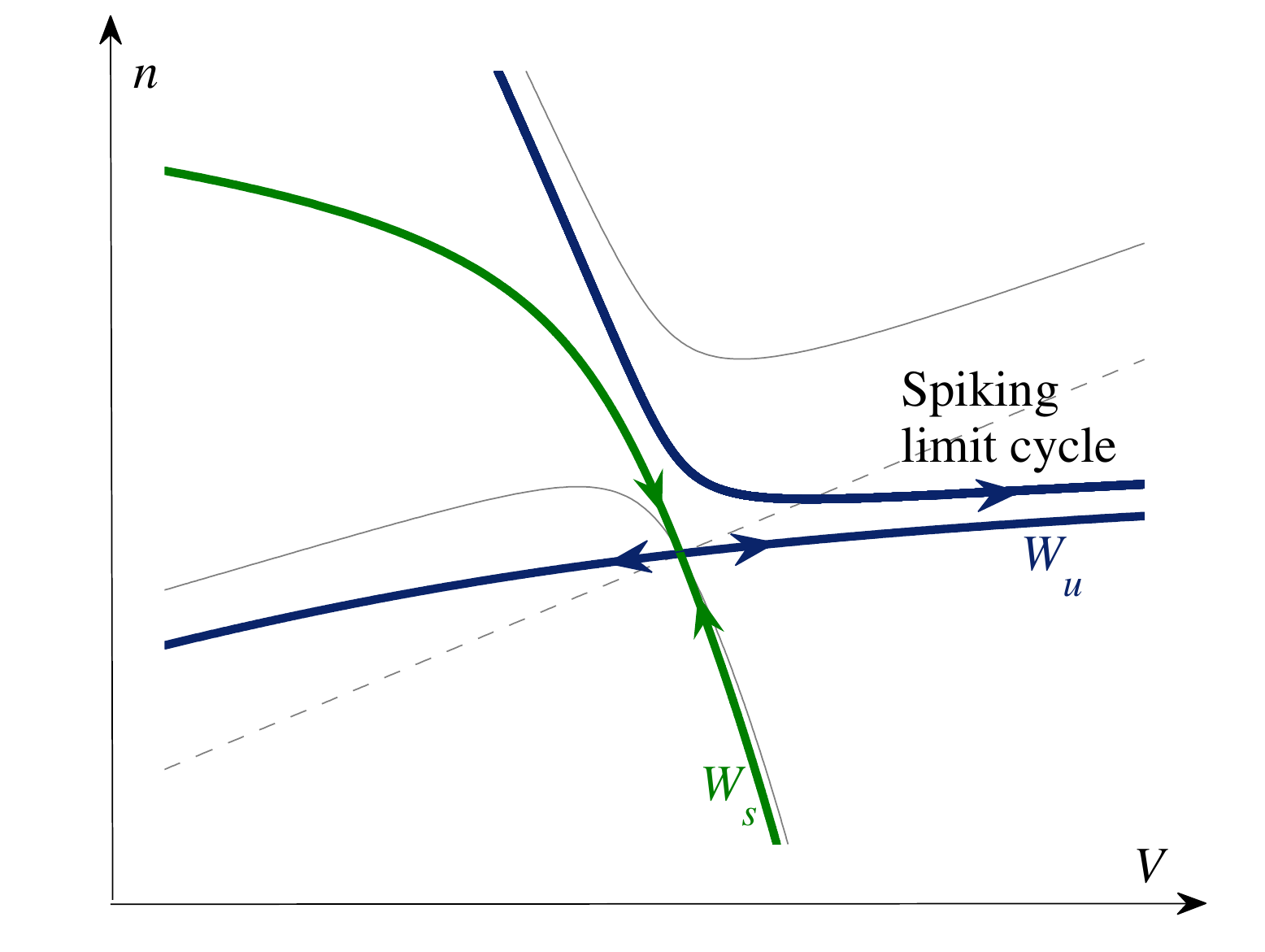}
}
\caption{{\bf Stable and unstable manifolds of the saddle point as $I_{app}$ varies.} (a,b) For $I_{app}<I_{SH}$, the left branch of the unstable manifold of the saddle $W_u$ (in violet) returns (after a spike generation) in the vicinity of the saddle point on the left of the stable manifold (in green) and toward rest. (c) At $I_{app}=I_{SH}$, there is a homoclinic orbit $\Gamma$ (in red), connecting the stable and the unstable manifolds. (d) For $I_{app}>I_{SH}$, the left branch of the unstable manifold returns in the vicinity of the saddle on the right of the stable manifold and converges toward the spiking limit cycle attractor.}\label{FIG: zoom saddle-homoclinic}
\end{figure}


Figure S\ref{FIG: zoom saddle-homoclinic} provides a close look of the saddle-homoclinic bifurcation in the reduced HH model (\ref{EQ: redHHCA}) with calcium current.
In order to rigorously prove the existence of this bifurcation, we rely on the normal form introduced in Lemma \ref{LEM: trans normal form} and exploit the timescale separation between $v$ and $w$ through geometrical singular perturbations theory (GSPT). To this aim we must briefly recall some basic of GSPT, relying on (S\ref{EQ: trans normal form}) as an explicit example. The interested reader will find in \cite{supp}{JONES95} an excellent introduction to the topic, and in \cite{supp}{KRSZ2001relax,KRSZ01,KRSZ01a} some recent extensions on which we rely for the forthcoming analysis.

\begin{figure}
\centering
\includegraphics[width=0.9\textwidth]{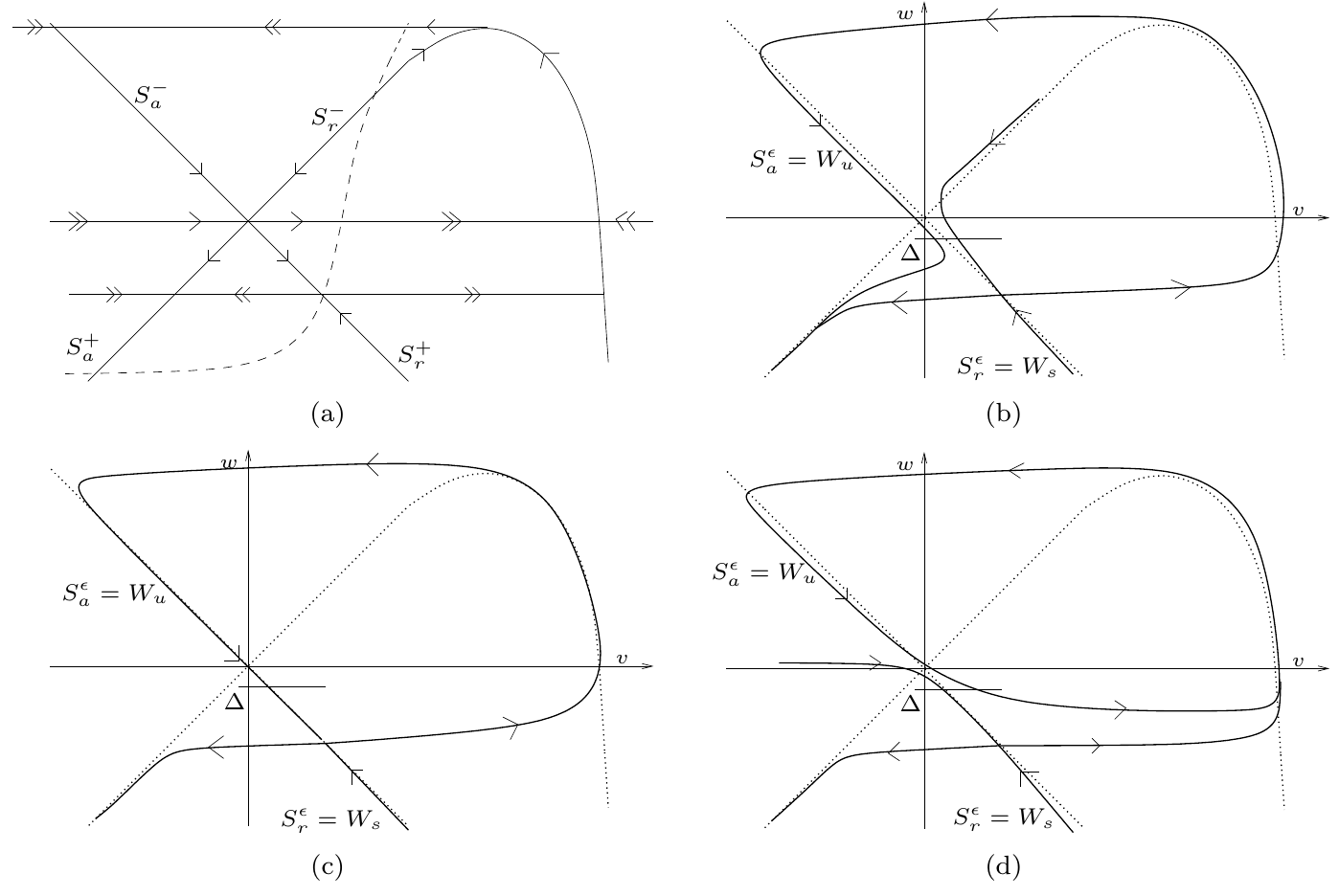}
\caption{{\bf Phase-portrait of (S\ref{EQ: trans normal form}) for different values of $\epsilon$ and $I$.} (a):  fast-slow dynamics of (S\ref{EQ: trans normal form}) for $\epsilon=I=0$. The attractive (resp. repelling) branches of the critical manifold $S_0$ above/below the transcritical singularity are denoted by $S_a^+/S_a^-$ (resp. $S_r^+/S_r^-$). (b): Continuation of the slow attractive $S_{a}^\epsilon$ and repelling $S_r^\epsilon$ manifolds for $\epsilon>0$ and $I<I_c(\sqrt{\epsilon})$, where $I_c(\sqrt{\epsilon})$ is defined as in Theorem \ref{THM: KRSZ01 thm reformulation}. (c): Continuation of $S_{a}^\epsilon$ and $S_r^\epsilon$ for $\epsilon>0$ and $I=I_c(\sqrt{\epsilon})$. (d): Continuation of $S_{a}^\epsilon$ and $S_r^\epsilon$ for $\epsilon>0$ and $I>I_c(\sqrt{\epsilon})$.}\label{FIG: GSPT full}
\end{figure}

The time rescaling $\tau:=\epsilon t$ transforms (S\ref{EQ: trans normal form}) into the equivalent system
\begin{subequations}\label{EQ: slow time scale}
\begin{eqnarray}
\epsilon\dot v&=&v^2-w^2+I+h_1(v,w,\epsilon)\IEEEyessubnumber\\
\dot w&=&(-1+h_2(v,w,\epsilon))\IEEEyessubnumber,
\end{eqnarray}
\end{subequations}
which describes the dynamics (S\ref{EQ: trans normal form}) in the slow timescale $\tau$.
In the limit $\epsilon=0$, commonly referred to as the \emph{singular limit}, one obtains from (S\ref{EQ: trans normal form}) and (S\ref{EQ: slow time scale}) two new dynamical systems: the \emph{reduced dynamics}
\begin{subequations}\label{EQ: reduced problem}
\begin{eqnarray}
0&=&v^2-w^2+I+h_1(v,w,\epsilon)\IEEEyessubnumber\\
\dot w&=&(-1+h_2(v,w,\epsilon))\IEEEyessubnumber,
\end{eqnarray}
\end{subequations}
evolves in the slow timescale $\tau$, while the \emph{layer dynamics}
\begin{subequations}\label{EQ: layer problem}
\begin{eqnarray}
\dot v&=&v^2-w^2+I+h_1(v,w,\epsilon)\IEEEyessubnumber\\
\dot w&=&0\IEEEyessubnumber,
\end{eqnarray}
\end{subequations}
evolves in the fast timescale $t$. Figure S\ref{FIG: GSPT full}(a) depicts the fast-slow dynamics associated to (S\ref{EQ: reduced problem})-(S\ref{EQ: layer problem}). The main idea behind GSPT is to combine the analysis of the reduced and layer dynamics to derive conclusion on the behavior of the nominal system, {\it i.e.} with $\epsilon>0$.

The reduced dynamics (S\ref{EQ: reduced problem}) is a dynamical system on the set $S_0:=\{(v,w)\in\RR^2:\ v^2-w^2+I+h_1(v,w,\epsilon)=0\}$, usually called the \emph{critical manifold}. The points in $S_0$ are indeed critical points of the layer dynamics (S\ref{EQ: layer problem}). More precisely, portions of $S_0$ on which $\frac{\partial\dot V}{\partial V}$ is non-vanishing are normally hyperbolic invariant manifolds of equilibria of the layer dynamics, whose stability is determined by the sign of $\frac{\partial\dot V}{\partial V}$. Conversely, points in $S_0$ where $\frac{\partial\dot V}{\partial V}=0$ constitute degenerate equilibria. In particular, the layer dynamics (S\ref{EQ: layer problem}) exhibits, for $I=0$, two degenerate equilibria. As depicted in Figure S\ref{FIG: GSPT full}(a), they are given by the self-intersection of the $V$-nullcline, corresponding to a transcritical singularity (Section \ref{SEC: trans bif}), and by the fold singularity at the maximum of the upper branch of the $V$-nullcline.

The basic result of GSPT, due to Fennichel \cite{supp}{FENICHEL79}, is that, for $\epsilon$ sufficiently small, non-degenerate portions of $S_0$ persist as nearby normally hyperbolic locally invariant manifolds $S^\epsilon$ of (S\ref{EQ: trans dyn syst}). More precisely, the \emph{slow manifold} $S^\epsilon$ lies in a neighborhood of $S_0$ of radius $O(\epsilon)$. The dynamics on $S^\epsilon$ is a small perturbation of the reduced dynamics (S\ref{EQ: reduced problem}). We point out that $S^\epsilon$ may not be unique, but is determined only up to $O(e^{-c/\epsilon})$, for some $c>0$. That is, two different choices of $S^\epsilon$ are exponentially close (in $\epsilon$) one to the other. Since the presented results are independent of the particular $S^\epsilon$ considered, we let this choice be arbitrary. The trajectories of the layer dynamics perturb to a stable and an unstable invariant foliations with basis $S^\epsilon$.

The analysis near degenerate points is more delicate. Only recently some works have treated this problem in its full generality for different types of singularities \cite{supp}{KRSZ2001relax,KRSZ01,KRSZ01a}. Figure S\ref{FIG: GSPT full} (b),(c),(d) sketch the extension of the attractive slow manifold $S_a^\epsilon$ after the fold point, and the three possible ways in which $S_a^\epsilon$ and the repelling slow manifold $S_r^\epsilon$ can continue after the trascritical singularity, depending on the injected current.

The result depicted in Figure S\ref{FIG: GSPT full} relies on the following theorem, adapted from \cite{supp}{KRSZ01}.

Let $\Delta:=\{(v,w)\in\mathbb R^2:\ v_-\leq v \leq v_+,\ w=\rho\}$, be the section depicted in Figure S\ref{FIG: GSPT full}, where $\rho<0$ and $|\rho|$ is sufficiently small, and $v_-,v_+$ are such that $\Delta\cap S_r^-\neq\emptyset$.
For a given $\epsilon>0$, let $q_{a,\epsilon}:=\Delta\cap S_a^\epsilon$ and $q_{r,\epsilon}:=\Delta\cap S^\epsilon_r$ be the intersections, whenever they exist, of respectively the attractive and repelling invariant submanifolds $S_a^\epsilon$ and $S_r^\epsilon$ with the section $\Delta$.
The following theorem reformulates in a compact way the discussion contained in Remark 2.2 and Section 3 of \cite{supp}{KRSZ01}\footnote{The first author is thankful to Prof. Szmolyan for his useful comments.} for systems with inputs of the form (S\ref{EQ: trans normal form}).
\begin{theorem}[Adapted from \cite{supp}{KRSZ01}]\label{THM: KRSZ01 thm reformulation}
Consider the system (S\ref{EQ: trans normal form}). Then there exists $\epsilon_0>0$ and a smooth function $I_c(\sqrt{\epsilon})$, defined on $[0,\epsilon_0]$ and satisfying $I_c(0)=0$, such that, for all $\epsilon\in(0,\epsilon_0]$, the following assertions hold
\begin{enumerate}
\item $q_{a,\epsilon}=q_{r,\epsilon}$ if and only if $I=I_c(\sqrt{\epsilon})$
\item there exists an open interval $A\ni I_c(\sqrt{\epsilon})$, such that, for all $I\in A$, it holds that $\Delta\cap S_a^\epsilon\neq\emptyset$, $\Delta\cap S_r^\epsilon\neq\emptyset$, and
$$\frac{\partial}{\partial I}\left(q_{a,\epsilon}-q_{r,\epsilon} \right)>0. $$
\end{enumerate}
\end{theorem}
Figure S\ref{FIG: GSPT full} illustrates this result.
\begin{remark}
The function $I_c(\sqrt{\epsilon})$ is related to the function $\lambda_c(\sqrt{\epsilon})$ defined in \cite[Remark 2.2]{supp}{KRSZ01} by $I_c(\sqrt{\epsilon}):=\epsilon\lambda_c(\sqrt{\epsilon})$. Similarly, given $\epsilon>0$, the parameter $I$ appearing in Theorem \ref{THM: KRSZ01 thm reformulation} is just the re-scaling $I=\epsilon\lambda$ of the parameter $\lambda$ appearing in \cite[Remark 2.2 and Sections 3]{supp}{KRSZ01}.
\end{remark}

Theorem \ref{THM: KRSZ01 thm reformulation} implies the existence of the saddle-homoclinic bifurcation in the reduced Hodgkin-Huxley model (\ref{EQ: redHHCA}) with calcium current. As stressed by Figure S\ref{FIG: GSPT full}(b,c,d), the slow attractive $S_a^\epsilon$ (resp. repelling $S_r^\epsilon$) manifold coincides with the unstable $W_u$ (resp. stable $W_s$) manifold of the saddle point, as it can be proved via qualitative arguments. Thus, for $I<I_c(\sqrt{\epsilon})$, the unstable manifold $W_u$ continues after the transcritical singularity on the left of $W_s$, toward the stable node. See Figure S\ref{FIG: zoom saddle-homoclinic}(a,b) and Figure S\ref{FIG: GSPT full}(b). For $I=I_c(\sqrt{\epsilon})$, $W_u$ extends after the transcritical point to $W_s$, forming the saddle-homoclinic trajectory, as depicted in Figure S\ref{FIG: zoom saddle-homoclinic}(c) and Figure S\ref{FIG: GSPT full}(c). For $I>I_c(\sqrt{\epsilon})$, the unstable manifold of the saddle $W_u$ continues after the transcritical singularity on the right of $W_s$, and spirals toward an exponentially stable limit cycle, whose existence can be proved with similar GSPT arguments (see for instance \cite{supp}{KRSZ2001relax}). This situation is the one depicted in Figure S\ref{FIG: zoom saddle-homoclinic}(d) and Figure S\ref{FIG: GSPT full}(d). 

\section{Hybrid singularly perturbed saddle-homoclinic bifurcation}
\label{SEC: hybrid SH bif}

The saddle-homoclinic bifurcation analysis provided in Section \ref{SEC: SP homoclinic bifurcation} for the reduced calcium-gated HH model naturally extends to the hybrid dynamics (\ref{EQ: IZHw2 model}) with $w_0<0$. Indeed, by construction, if $w_0<0$, this model can be transformed in the normal form (S\ref{EQ: trans normal form}) derived in Lemma \ref{LEM: trans normal form}, and Theorem \ref{THM: KRSZ01 thm reformulation} applies directly.

\begin{figure}
\centering
\includegraphics[width=0.9\textwidth]{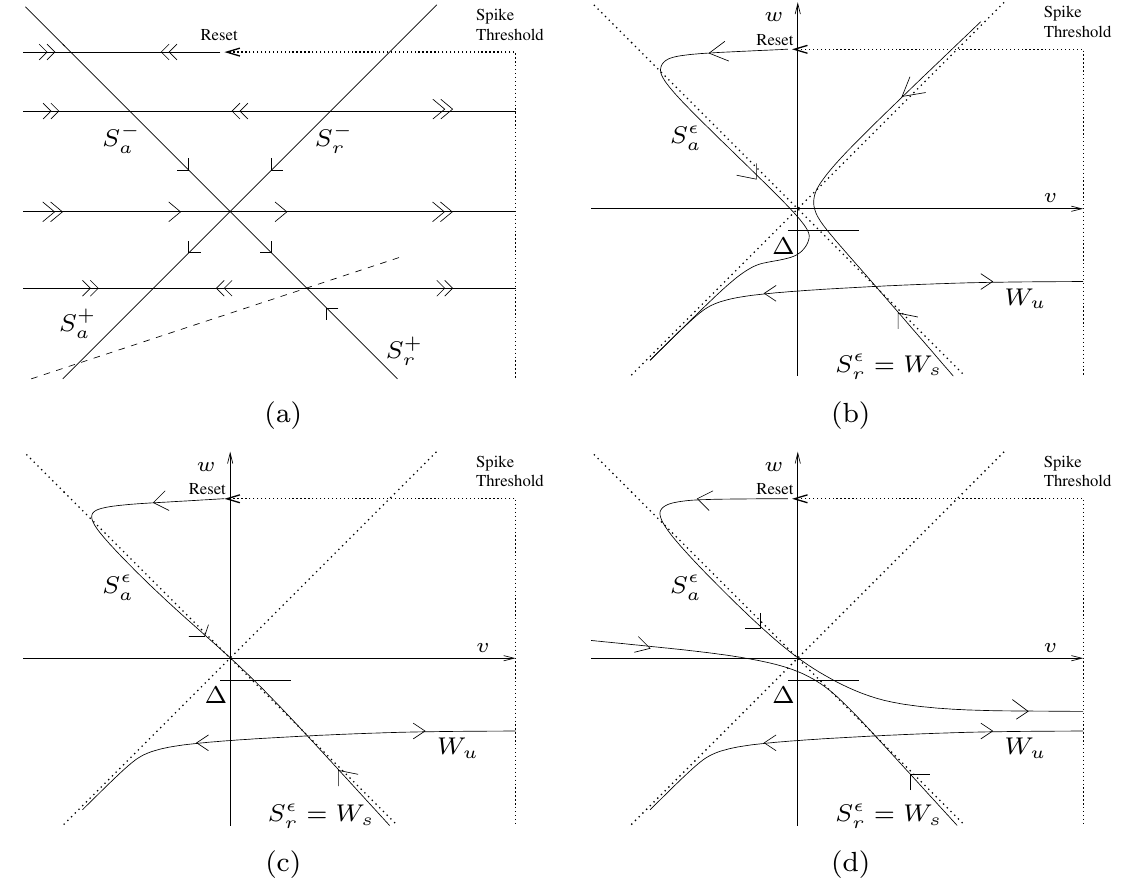}
\caption{{\bf Phase-portrait of the hybrid dynamics (\ref{EQ: IZHw2 model}) for different values of $\epsilon$ and $I$.} (a): Fast-slow dynamics for $\epsilon=I=0$. The attractive (resp. repelling) branches of the critical manifold $S_0$ above/below the transcritical singularity are denoted by $S_a^+/S_a^-$ (resp. $S_r^+/S_r^-$). (b): Continuation of the slow manifolds for $\epsilon>0$ and $I<I_c(\sqrt{\epsilon})$, where $I_c(\sqrt{\epsilon})$ is defined as in Theorem \ref{THM: KRSZ01 thm reformulation}. (c): Continuation of the slow manifolds for $\epsilon>0$ and $I=I_c(\sqrt{\epsilon})$. (d): Continuation of the slow manifolds for $\epsilon>0$ and $I>I_c(\sqrt{\epsilon})$. }\label{FIG: hybrid GSPT}
\end{figure}

Similarly to the derivation in Section \ref{SEC: SP homoclinic bifurcation}, we can associate to (\ref{EQ: IZHw2 model}) two new dynamical systems, describing its singular dynamics in the fast timescale $t$ and in the slow timescale $\tau:=\epsilon t$, respectively. More precisely, the \emph{hybrid reduced dynamics}
\begin{IEEEeqnarray*}{rClCl}
0&=&v^2-w^2+I&\quad\quad\quad&\text{if}\ v\geq v_{th},\ \text{then}\\
\dot w&=&av-w+w_0&\quad\quad\quad&v\leftarrow c,\ u\leftarrow d,
\end{IEEEeqnarray*}
evolves in the slow time scale $\tau$, while the \emph{hybrid layer dynamics}
\begin{IEEEeqnarray*}{rClCl}
\dot v&=&v^2-w^2+I&\quad\quad\quad&\text{if}\ v\geq v_{th},\ \text{then}\\
\dot w&=&0&\quad\quad\quad&v\leftarrow c,\ u\leftarrow d,
\end{IEEEeqnarray*}
evolves in the fast time-scale $t$. Figure S\ref{FIG: hybrid GSPT}(a) depicts the associated slow-fast hybrid dynamics for $I=0$, $a\in(0,1)$, and $w_0<0$.

The analysis of the non-singular limit follows the same line as the analysis developed in Section \ref{SEC: SP homoclinic bifurcation} for the continuous time case. The only difference is that the return mechanism, provided in the continuous time case by the right attractive branch of the critical manifold $S_0$, is now replaced by the hybrid reset. The result is summarized in Figure S\ref{FIG: hybrid GSPT} (b),(c),(d). The slow attractive manifold $S_a^\epsilon$ is chosen as the continuation of the trajectory starting at the reset point. In this way, it also coincides with the image of the unstable manifold of the saddle $W_u$ through the hybrid reset map. As in the continuous time case, the stable manifold of the saddle $W_s$ can be shown to coincide with the slow repelling manifold $S_r^\epsilon$. Let $I_c(\sqrt{\epsilon})$ be defined as in Theorem \ref{THM: KRSZ01 thm reformulation}. For $I<I_c(\sqrt{\epsilon})$, the unstable manifold $W_u$ is brought back on the left of $W_s$, and the resting state is globally stable, as in Figure S\ref{FIG: hybrid GSPT}(b). At $I=I_c(\sqrt{\epsilon})$, the hybrid reset connects $W_u$ and $W_s$, corresponding to a hybrid homoclinic bifurcation. The associated phase-portrait is shown in Figure S\ref{FIG: hybrid GSPT}(c). Finally, for $I>I_c(\sqrt{\epsilon})$, the unstable manifold $W_u$ is brought by the reset mechanism on the right of $W_s$ and directly into the newborn hybrid limit cycle attractor, as in Figure S\ref{FIG: hybrid GSPT}(d).

We emphasize that the existence of the hybrid homoclinic bifurcation and the associated critical value $I_c$ are independent of the reset point $(c,d)$, provided that the trajectory starting from $(c,d)$ be attracted toward the left upper branch of the critical manifold, as in Figure S\ref{FIG: hybrid GSPT}. As discussed in Section \ref{SEC: SP homoclinic bifurcation}, under this condition, two different choices of the reset point are associated to two slow manifolds $S_a^\epsilon$ that are $O(e^{-c/\epsilon})$ near, for some $c>0$. By the result in Theorem \ref{THM: KRSZ01 thm reformulation}, the two values of the critical current $I_c$ for which the slow attractive manifold extends to the slow repelling manifold are again $O(e^{-c/\epsilon})$ near. This also ensures that the value for which the hybrid homoclinic bifurcation happens is independent of the reset point, modulo variations that are $O(e^{-c/\epsilon})$. The same robustness properties are not shared by fold hybrid models (see Section \ref{SEC: robust ADPs} below).

\section{Phase-plane analysis of the transcritical hybrid TC neuron model}
\label{SEC: hybrid TC phase-plane}

\begin{figure}[htbp]
\begin{center}
\includegraphics[width=0.9\textwidth]{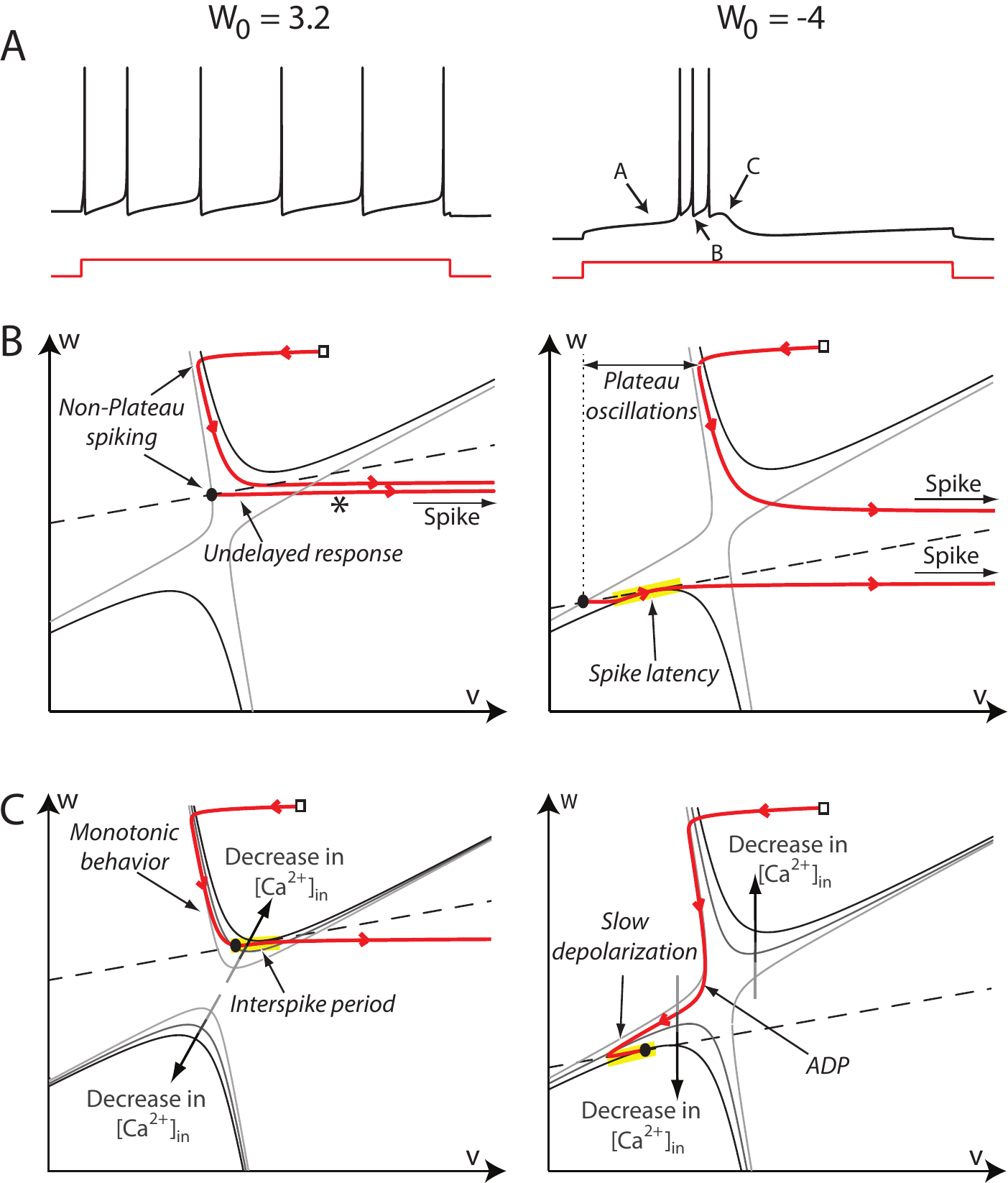}
\end{center}
\caption{{\bf Transcritical hybrid modeling of TC cell (\ref{EQ: quantitative IZHw2 model}) in the low ($w_0 = 3.2$, left) and high calcium conductance modes ($w_0 = -4$, right).} (\textbf{a}) Time-courses of the model in both conditions. ({\bf b} and {\bf c}) Phase-portrait of the transcritical hybrid model in both conditions. Trajectories are depicted as solid oriented red lines. The reset point is depicted as a square $\Box$, while the (instantaneous) hyperpolarized state as a filled circle $\bullet$. The $w$-nullcline is depicted as a dashed line. In (\textbf{b}), the gray (black) thin solid line is the $v$-nullcline when the current step in off (on). In (\textbf{c}), the $v$-nullcline is depicted as gray thin lines of different darkness. As sketched in the figure, light gray correspond to large calcium concentrations, whereas dark gray to small.}
\label{fig:s4}
\end{figure}

Fig. S\ref{fig:s4} shows the time-course of the qualitative transcritical hybrid modeling of thalamocortical relay (TC) neurons (\ref{EQ: quantitative IZHw2 model}) in the low and high calcium conductance modes, as well as the corresponding phase portraits. The analysis of the phase-portraits gives insights on the mechanisms by which TC cells exhibit tonic firing or bursting when submitted to a similar step input, according to the initial resting potential.

When $w_0>0$ (Fig. S\ref{fig:s4}, left), which corresponds to a small calcium conductance (T-type calcium channels are inactivated), the hyperpolarized state belongs to the upper branch of the $v$-nullcline. Application of a depolarizing current step lifts the voltage nullcline above the resting state, thus generating a transient non-delayed action potential (marked with a $*$ in Fig. S\ref{fig:s4}{\bf b}, left). Note that, when the hyperpolarized state belongs to the upper branch of the $v$-nullcline, no plateau oscillations are possible (Fig. S\ref{fig:s4}{\bf b}, left). Furthermore, the relaxation toward the hyperpolarized state is necessarily monotone ({\it i.e.} no ADP), as stressed in Fig. S\ref{fig:s4}{\bf c}, left.

At the generation of the first spike, a certain amount of calcium ions enter the cell due to the presence of high voltage activated calcium currents. The subsequent activation of outward calcium pump currents\footnote{This activation is modeled by the hybrid reset $z\leftarrow z+d_z$.} transiently hyperpolarizes the cell. As the intracellular calcium is expelled (Fig. \ref{fig:s4}{\bf c}, left), calcium pump currents slowly deactivate and the cell slowly depolarizes (interspike period), until the spiking threshold is reached and a new action potential is fired.

When $w_0<0$ (Fig. S\ref{fig:s4}, right), which corresponds to a large calcium conductance (T-type calcium channels are available), the hyperpolarized state belongs to the lower branch of the $v$-nullcline, and is more hyperpolarized than for $w_0>0$, as in experiments. When a depolarizing current step is applied, this branch falls below the $w$-nullcline (Fig. S\ref{fig:s4}{\bf b}, right). In order to generate the first spike, the state travels in the narrow region between the two nullclines, resulting in a pronounced latency. Furthermore, the relative position of the hyperpolarized state (lower branch) with respect to (hybrid) spiking cycle (upper branch) clearly explains the generation mechanism of plateau oscillations. These high frequency plateau oscillations (burst) continue until the intracellular calcium, which enters at each action potential, and the associated pumps current ({\it i.e.} $z$) are sufficiently large. Plateau oscillations then terminate in a (hybrid) saddle-homoclinic bifurcation (Fig. S\ref{fig:s4}{\bf c}, right). At the end of the burst, the system converges toward the hyperpolarized state following the left branch of the $v$-nullcline, thus generating a marked ADP at the passage near the nullcline's funnel. The subsequent slow phase is mainly ruled by the variations of the intracellular calcium. With the adopted simple dynamics it consists in a slow depolarization that follows the decrease of the intracellular calcium (Fig. S\ref{fig:s4}{\bf c}, right). A finer and more physiological modeling of the intracellular calcium dynamics could reproduce {\it in vitro} recordings with a higher degree of fidelity.

\section{Modeling TC neurons with a fold hybrid model}
\label{EQ: fold hybrid model}

\begin{figure}[htbp]
\begin{center}
	\includegraphics[width=0.6\textwidth]{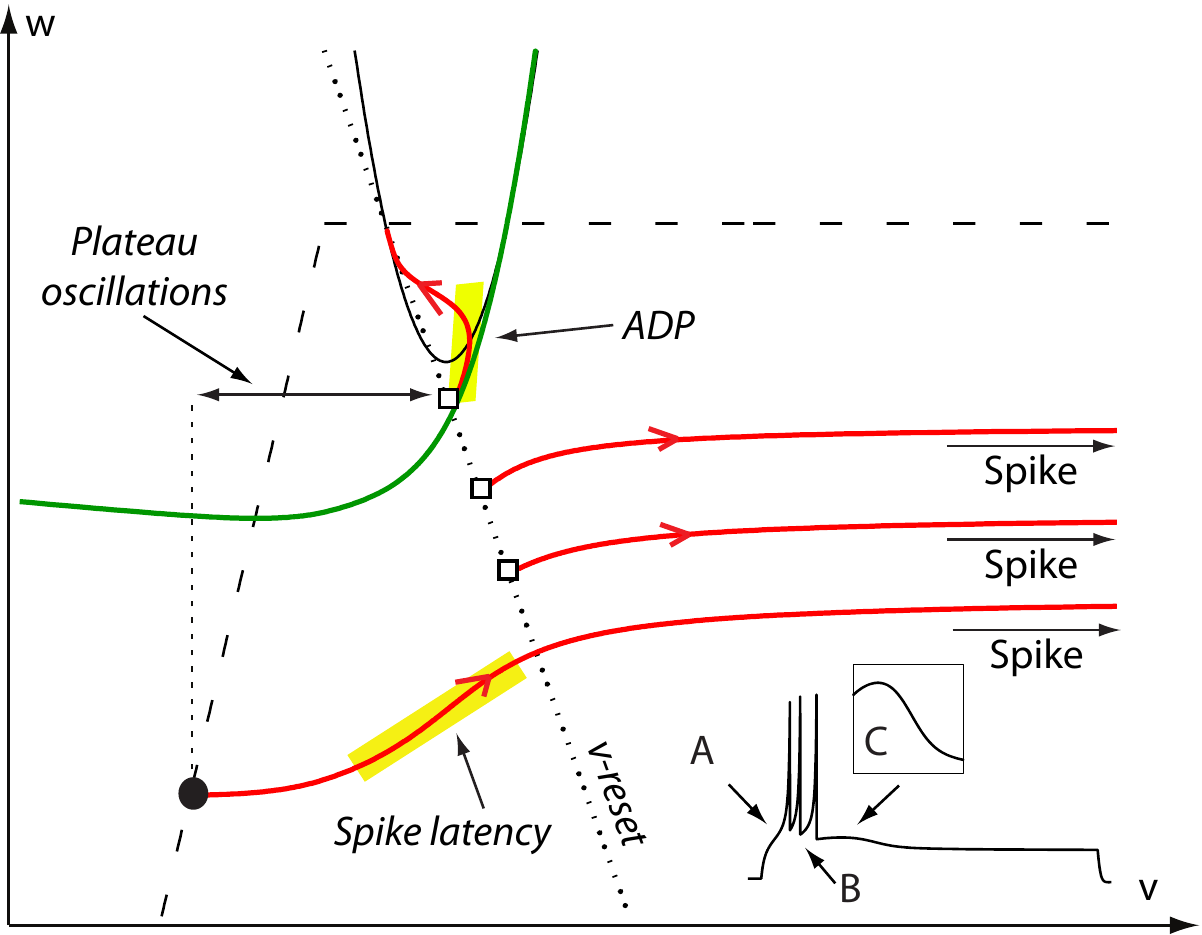}
\end{center}
\caption{{\bf Phase portrait and voltage time course (inset) of the Izhikevich model of TC neurons.} The trajectory is depicted as a solid oriented red line. The $v$-nulcline is depicted a solid line, the $w$-nullcline as a dashed line. The point line sketches the $v$-rest law. The stable manifold of the saddle point is depicted as a green line. Parameters as in \cite[Sections 8.3.1]{supp}{IZHIKEVICH2007}.}\label{FIG: izhi TC_like}
\end{figure}

This section reproduces the modeling of TC neuron with the fold hybrid model. Such a model was recently proposed by Izhikevich and used in a large-scale model of the mammalian thalamic system \cite{supp}{izhikevich2008large}. The model succeeds in reproducing the firing pattern of TC neurons by modifying the reset mechanism. Fig. S\ref{FIG: izhi TC_like} shows that the voltage time-course reproduces with fidelity the three hallmarks of large calcium conductances, that is (A) spike latency, (B) plateau oscillations, and (C) ADP. However, the phase portrait illustrates that the modified $v$-reset law has no obvious physiological interpretation. As a consequence, plateau oscillations are not endogenously generated in this model. Furthermore, as opposed to both the quantitative model and the proposed hybrid model, the ADP is generated following an increase of the recovery variable, when the reset point crosses the stable manifold of the saddle.

\section{Robust ADP generation}
\label{SEC: robust ADPs}

\begin{figure}[htbp]
\begin{center}
	\includegraphics[width=0.8\textwidth]{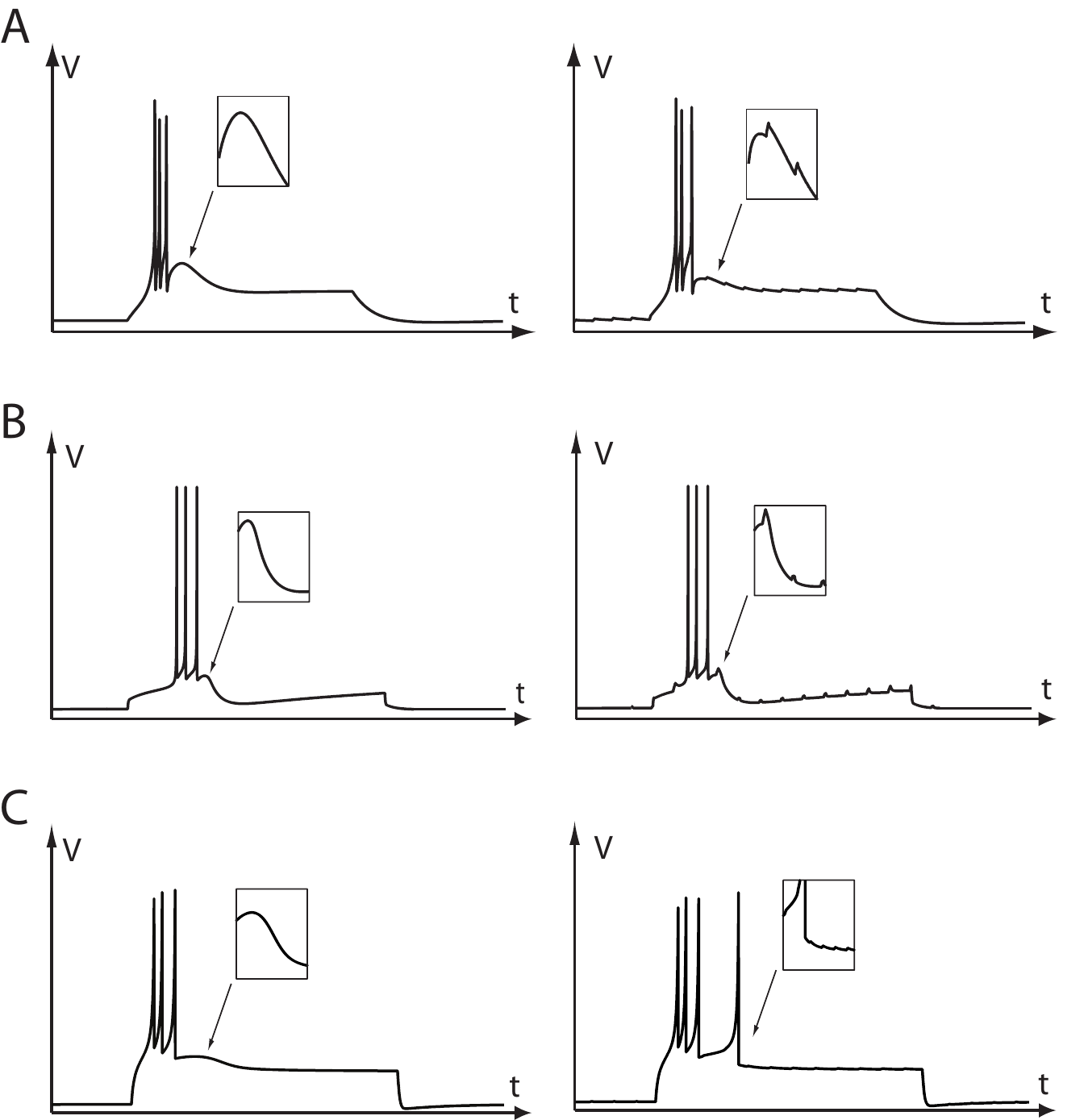}
\end{center}
\caption{{\bf Nominal step response (left) and step response in the presence of small current pulses in the 200 compartments TC neuron model ({\bf a}), the proposed hybrid model  ({\bf b}), and the Izhikevich model of TC neuron (\cite[Section 8.3.1]{supp}{IZHIKEVICH2010}) ({\bf c}).}}
\label{FIG: TC compartment}\label{FIG: pacemaking comparison}
\end{figure}

For a neuron model to be biologically relevant, it should be robust to exogenous disturbances (small synaptic inputs, thermal noise, etc.). The firing pattern, in particular, should remain unchanged. Figure S\ref{FIG: TC compartment} compares the perturbation robustness of three TC neuron models to small current impulses. It suggests that the fold hybrid model is less robust than the transcritical model, because a tiny pulse is sufficient to generate an extra action potential at the ADP apex.

\begin{figure}[htbp]
\begin{center}
	\includegraphics[width=0.8\textwidth]{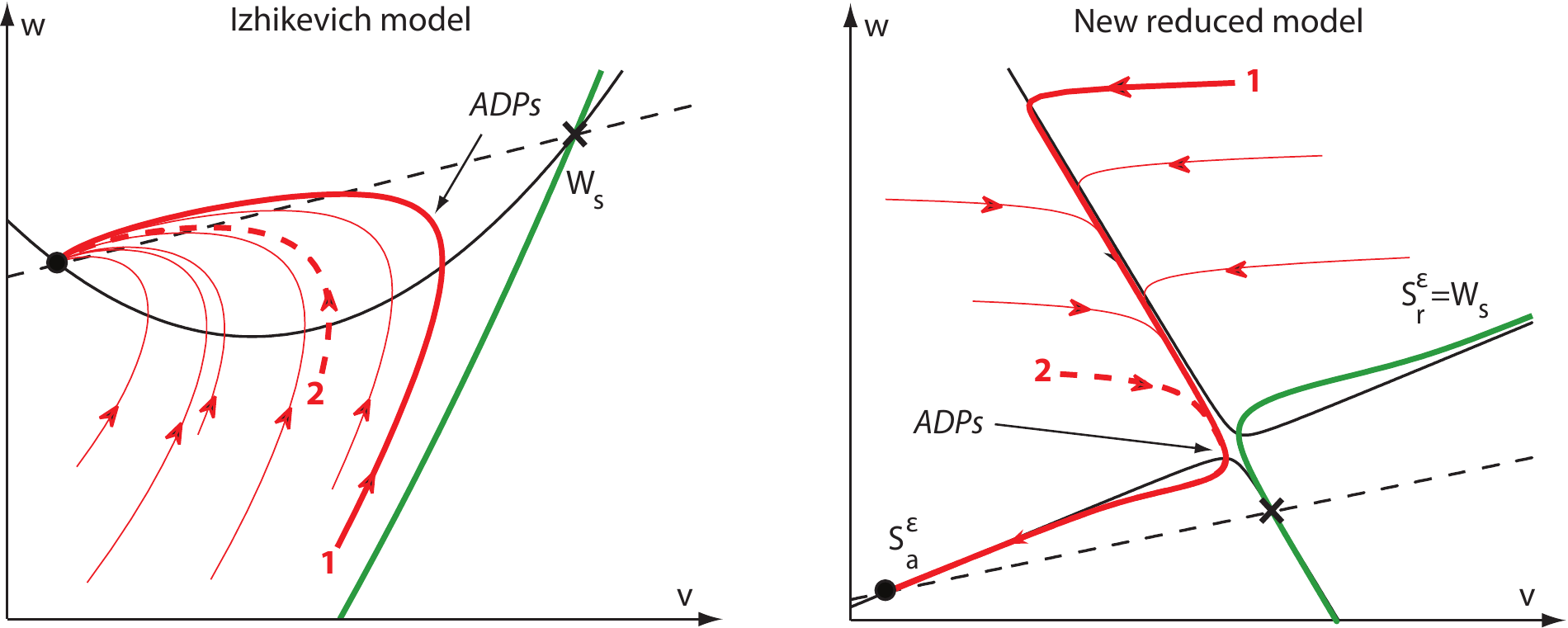}
\end{center}
\caption{{\bf Comparison of the ADP generation mechanisms in the fold (left) and in the transcritical hybrid models (right).} The stable manifold of the saddle ($\times$) is depicted in green. In the Izhikevich model ADPs are generated by sliding near the stable manifold of the saddle and crossing the $v$-nullcline from below. In our proposed hybrid model, ADP are robustly generated along the attractor $S_a^\epsilon$.}\label{FIG: ADP generation comparison}
\end{figure}

The difference in robustness is explained by the different ADP generation mechanisms, as illustrated in Figure S\ref{FIG: ADP generation comparison}. In the fold model, ADPs are generated when trajectories cross the $v$-nullcline from below (cf. Figures S\ref{FIG: izhi TC_like} and S\ref{FIG: ADP generation comparison}). The absence of any robust attractor in the ADP generation region makes the ADP height and shape heavily dependent on the exact reset point. Moreover, when small current pulses are applied, the ADP generation is disrupted, and the model fires an extra (non-physiological) spike.

Conversely, ADP generation in the transcritical hybrid model is robustly governed by the attractor $S_a^\epsilon$ that steers the trajectories through the ADP apex and toward the resting point. That is the reason why the ADP height and shape barely depend on chosen reset point. At the same time, the persistence to small perturbations of this invariant manifold \cite{supp}{HIPUSH77} ensures, as required in biologically meaningful conditions, the robustness of the ADP generation mechanism to small inputs.

\clearpage

\bibliographystyle{supp}{plain}
\bibliography{supp}{../../../refs}{Supplementary references}

\newcommand{\SortNoop}[1]{} 
  \def\nesic{Ne\v{s}i\'{c}\,} %
  \def\astrom{{\SortNoop{As}\AA}str{\"{o}}m\,}\let\c=\cedille
\begin{thebibliography}{10}

\bibitem{Azouz01041996}
R.~Azouz, M.~S. Jensen, and Y.~Yaari.
\newblock Ionic basis of spike after-depolarization and burst generation in
  adult rat hippocampal {CA1} pyramidal cells.
\newblock {\em The Journal of Physiology}, 492(1):211--223, 1996.

\bibitem{Beurrier15011999}
C.~Beurrier, P.~Congar, B.~Bioulac, and C.~Hammond.
\newblock Subthalamic nucleus neurons switch from single-spike activity to
  burst-firing mode.
\newblock {\em The Journal of Neuroscience}, 19(2):599--609, 1999.

\bibitem{canavier2006increase}
C.C. Canavier and R.S. Landry.
\newblock An increase in {AMPA} and a decrease in {SK} conductance increase
  burst firing by different mechanisms in a model of a dopamine neuron in vivo.
\newblock {\em Journal of neurophysiology}, 96(5):2549--2563, 2006.

\bibitem{catterall2005international}
W.~A. Catterall, E.~Perez-Reyes, T.~P. Snutch, and J.~Striessnig.
\newblock {International Union of Pharmacology. XLVIII}. {Nomenclature} and
  structure-function relationships of voltage-gated calcium channels.
\newblock {\em Pharmacological reviews}, 57(4):411, 2005.

\bibitem{TJP:TJP2775}
S.~Chen and Y.~Yaari.
\newblock Spike {Ca2+} influx upmodulates the spike afterdepolarization and
  bursting via intracellular inhibition of {KV7/M} channels.
\newblock {\em The Journal of Physiology}, 586(5):1351--1363, 2008.

\bibitem{destexhe1998dendritic}
A.~Destexhe, M.~Neubig, D.~Ulrich, and J.~Huguenard.
\newblock Dendritic low-threshold calcium currents in thalamic relay cells.
\newblock {\em The Journal of neuroscience}, 18(10):3574--3588, 1998.

\bibitem{DRMASESE11}
G.~Drion, L.~Massotte, R.~Sepulchre, and V.~Seutin.
\newblock How modeling can reconcile apparently discrepant experimental
  results: The case of pacemaking in dopaminergic neurons.
\newblock {\em PLoS Comput Biol}, 7(5):e1002050, 05 2011.

\bibitem{ERTE10}
G.~B. Ermentrout and D.~H. Terman.
\newblock {\em Mathematical Foundations of Neuroscience}.
\newblock Interdisciplinary Applied Mathematics. Springer, 2010.

\bibitem{fitzhugh61}
R.~FitzHugh.
\newblock Impulses and physiological states in theoretical models of nerve
  membrane.
\newblock {\em Biophysical J.}, 1:445--466, 1961.

\bibitem{halnes2011multi}
G.~Halnes, S.~Augustinaite, P.~Heggelund, G.~T. Einevoll, and M.~Migliore.
\newblock A multi-compartment model for interneurons in the dorsal lateral
  geniculate nucleus.
\newblock {\em BMC Neuroscience}, 12(Suppl 1):P222, 2011.

\bibitem{hille1991ionic}
B.~Hille.
\newblock {\em Ionic Channels of Excitable Membranes}.
\newblock Sinauer, Sunderland, Massachusetts, 2nd. edition, 1991.

\bibitem{HODHUX}
A.~Hodgkin and A.~Huxley.
\newblock A quantitative description of membrane current and its application to
  conduction and excitation in nerve.
\newblock {\em J. Physiol}, 117:500--544, 1952.

\bibitem{IZHIKEVICH03}
E.~M. Izhikevich.
\newblock A simple model of spiking neurons.
\newblock {\em IEEE Trans. on Neural Networks}, 14(6):1569--1572, 2003.

\bibitem{IZHIKEVICH2007}
E.~M. Izhikevich.
\newblock {\em Dynamical Systems in Neuroscience: The Geometry of Excitability
  and Bursting}.
\newblock MIT Press, 2007.

\bibitem{IZHIKEVICH2010}
E.~M. Izhikevich.
\newblock Hybrid spiking models.
\newblock {\em Phil. Trans. R. Soc. A}, 368:5061--5070, 2010.

\bibitem{izhikevich2008large}
E.M. Izhikevich and G.M. Edelman.
\newblock Large-scale model of mammalian thalamocortical systems.
\newblock {\em Proceedings of the national academy of sciences},
  105(9):3593--3598, 2008.

\bibitem{McCormick92a}
D~A McCormick and J~R Huguenard.
\newblock A model of the electrophysiological properties of thalamocortical
  relay neurons.
\newblock {\em J Neurophysiol}, 68(4):1384--1400, 1992.

\bibitem{Molineux23112005}
M.~L. Molineux, F.~R. Fernandez, W.~H. Mehaffey, and R.~W. Turner.
\newblock {A-Type and T-Type} currents interact to produce a novel spike
  latency-voltage relationship in cerebellar stellate cells.
\newblock {\em The Journal of Neuroscience}, 25(47):10863--10873, 2005.

\bibitem{pospischil2011comparison}
M.~Pospischil, Z.~Piwkowska, T.~Bal, and A.~Destexhe.
\newblock Comparison of different neuron models to conductance-based
  post-stimulus time histograms obtained in cortical pyramidal cells using
  dynamic-clamp in vitro.
\newblock {\em Biological cybernetics}, 105(2):167--180, 2011.

\bibitem{Rekling01111997}
J.~C. Rekling and J.~L. Feldman.
\newblock Calcium-dependent plateau potentials in rostral ambiguus neurons in
  the newborn mouse brain stem in vitro.
\newblock {\em Journal of Neurophysiology}, 78(5):2483--2492, 1997.

\bibitem{richert2011efficient}
M.~Richert, J.M. Nageswaran, N.~Dutt, and J.L. Krichmar.
\newblock An efficient simulation environment for modeling large-scale cortical
  processing.
\newblock {\em Frontiers in Neuroinformatics}, 5(19), 2011.

\bibitem{Rinzel:1989:ANE:94605.94613}
J.~Rinzel and G.~B. Ermentrout.
\newblock {\em Analysis of neural excitability and oscillations}, pages
  135--169.
\newblock MIT Press, Cambridge, MA, USA, 1989.

\bibitem{SEYDEL94}
R.~Seydel.
\newblock {\em Practical bifurcation and stability analysis}, volume~5 of {\em
  Interdisciplinary Applied Mathematics}.
\newblock Springer-Verlag, New York, third edition, 2010.

\bibitem{Sherman01a}
S~M Sherman.
\newblock Tonic and burst firing: dual modes of thalamocortical relay.
\newblock {\em Trends Neurosci}, 24(2):122--126, 2001.

\bibitem{wang1991model}
X.~J. Wang, J.~Rinzel, and M.~A. Rogawski.
\newblock A model of the {T}-type calcium current and the low-threshold spike
  in thalamic neurons.
\newblock {\em Journal of neurophysiology}, 66(3):839--850, 1991.

\end{thebibliography}



\begin{thebibliography}{10}

\bibitem{ERMEXPPbook}
G.~B. Ermentrout.
\newblock {\em Simulating, analyzing, and animating dynamical systems: a guide
  to XPPAUT for researchers and students}.
\newblock SIAM Press, Philadelphia, PA, USA, 2002.

\bibitem{FENICHEL79}
N.~Fenichel.
\newblock Geometric singular perturbation theory.
\newblock {\em J. Diff. Eq.}, 31:53--98, 1979.

\bibitem{HIPUSH77}
M.~Hirsch, C.~Pugh, and M.~Shub.
\newblock {\em Invariant Manifolds}.
\newblock Lecture Notes in Mathematics. Springer-Verlag, Berlin, Germany, 1977.

\bibitem{IZHIKEVICH2007}
E.~M. Izhikevich.
\newblock {\em Dynamical Systems in Neuroscience: The Geometry of Excitability
  and Bursting}.
\newblock MIT Press, 2007.

\bibitem{IZHIKEVICH2010}
E.~M. Izhikevich.
\newblock Hybrid spiking models.
\newblock {\em Phil. Trans. R. Soc. A}, 368:5061--5070, 2010.

\bibitem{izhikevich2008large}
E.M. Izhikevich and G.M. Edelman.
\newblock Large-scale model of mammalian thalamocortical systems.
\newblock {\em Proceedings of the national academy of sciences},
  105(9):3593--3598, 2008.

\bibitem{JONES95}
C.~K.~R. Jones.
\newblock Geometric singular perturbation theory.
\newblock In {\em Dynamical systems. Springer Lecture Notes in Math. 1609},
  pages 44--120, Berlin, 1995. Springer.

\bibitem{KRSZ01a}
M.~Krupa and P.~Szmolyan.
\newblock Extending geometrical singular perturbation theory to nonhyperbolic
  points - folds and canards points in two dimensions.
\newblock {\em SIAM J. Math. Analysis}, 33(2):286--314, 2001.

\bibitem{KRSZ01}
M.~Krupa and P.~Szmolyan.
\newblock Extending slow manifolds near transcritical and pitchfork
  singularities.
\newblock {\em Nonlinearity}, 14:1473--1491, 2001.

\bibitem{KRSZ2001relax}
M.~Krupa and P.~Szmolyan.
\newblock Relaxation oscillation and canard explosion.
\newblock {\em J. Differential Equations}, 174(2):312--368, 2001.

\bibitem{Lee_v2006}
J.~Lee.
\newblock {\em Introduction to smooth manifolds}.
\newblock Graduate Texts in Mathematics. Springer-Verlag, Berlin, Germany,
  2006.

\bibitem{SEYDEL94}
R.~Seydel.
\newblock {\em Practical bifurcation and stability analysis}, volume~5 of {\em
  Interdisciplinary Applied Mathematics}.
\newblock Springer-Verlag, New York, third edition, 2010.

\bibitem{strogatz00}
S.~H. Strogatz.
\newblock From {Kuramoto} to {Crawford}: exploring the onset of synchronization
  in population of coupled oscillators.
\newblock {\em Physica D}, 143:1--20, 2000.

\end{thebibliography}

\end{document}